 \documentclass[final,3p,times]{elsarticle}
\usepackage{hyperref}  
 
 \hypersetup{
     colorlinks=false,
     linkcolor=black,
     citecolor=black,
     filecolor=black,
     urlcolor=black,
 }
\usepackage{graphicx}
\usepackage{caption}
\usepackage{subcaption}

\usepackage{siunitx}

\usepackage{amsmath,amssymb,amsthm}

\usepackage[numbers]{natbib}
\let\cite=\citet

\newcommand{\bef}{{\bf f}}
\newcommand{\bg}{{\bf g}}

\newcommand{\bn}{{\bf n}}

\newcommand{\bp}{{\bf p}}

\newcommand{\bu}{{\bf u}}
\newcommand{\bv}{{\bf v}}
\newcommand{\bw}{{\bf w}}
\newcommand{\bx}{{\bf x}}
\newcommand{\by}{{\bf y}}
\newcommand{\bz}{{\bf z}}

\newcommand{\bF}{{\bf F}}

\newcommand{\bI}{{\bf I}}

\newcommand{\bP}{{\bf P}}

\newcommand{\bU}{{\bf U}}
\newcommand{\bV}{{\bf V}}
\newcommand{\bW}{{\bf W}}
\newcommand{\bX}{{\bf X}}

\newcommand{\bnu}{{\boldsymbol \nu}}

\newcommand{\btau}{{\boldsymbol \tau}}

\newcommand{\calC}{{\mathcal C}}

\newcommand{\polQ}{{\mathbb Q}}

\newtheorem{rmk}{Remark}

\usepackage{cases}
\usepackage{tikz, xcolor}
\usepackage{xspace}
\usetikzlibrary{shapes.geometric}
\usetikzlibrary{shapes,arrows}

\tikzstyle{block} = [rectangle, draw, text centered, rounded corners, 
                     text width=7em,    
                     minimum height=3em, node distance=3.3cm]
\tikzstyle{line} = [draw, -latex']
\tikzstyle{cloud} = [draw, ellipse, node distance=2.5cm, minimum height=3em]

\newcommand{\step}[1]{\noindent\raisebox{1.5pt}[10pt][0pt]{\tiny\framebox{$#1$}}\xspace}
\newcommand{\alt}[1]{}

\newcommand{\ie}{i.e.,\@\xspace}
\newcommand{\eg}{e.g.\@\xspace}

\newcommand{\Real}{\mathbb R}
\newcommand{\CROSS}{{\times}}   
\newcommand{\SCAL}{{\cdot}}       

\journal{Journal of Computational Physics}

\begin{document}

\begin{frontmatter}



\title{Numerical Simulations of Bouncing Jets}


\author[tamu]{Andrea Bonito\fnref{fnab}}
\ead{bonito@math.tamu.edu}
\author[tamu]{Jean-Luc Guermond}
\ead{guermond@math.tamu.edu}
\author[tamu]{Sanghyun Lee\fnref{fnab}}
\ead{shlee@math.tamu.edu}

\fntext[fnab]{Partially supported by the NSF grant DMS-1254618 and by award
KUS-C1-016-04 made by King Abdullah University of Science and Technology (KAUST)}

\address[tamu]{Department of Mathematics, Texas A\&M University,\\ College Station, Texas 77843-3368, USA}

\begin{abstract}
  Bouncing jets are fascinating phenomenons occurring under certain
  conditions when a jet impinges on a free surface. This effect is
  observed when the fluid is Newtonian and the jet falls in a bath
  undergoing a solid motion.  It occurs also for non-Newtonian fluids
  when the jets falls in a vessel at rest containing the same fluid.
  We investigate numerically the impact of the experimental setting
  and the rheological properties of the fluid on the onset of the
  bouncing phenomenon.  Our investigations show that the occurrence of
  a thin lubricating layer of air separating the jet and the rest of
  the liquid is a key factor for the bouncing of the jet to happen.
  The numerical technique that is used consists of a projection method
  for the Navier-Stokes system coupled with a level set formulation
  for the representation of the interface. The space approximation is
  done with adaptive finite elements. Adaptive refinement is shown to
  be very important to capture the thin layer of air that is
  responsible for the bouncing.
\end{abstract}


\end{frontmatter}


\section{Introduction} 
\label{sec:intro}

The ability of a jet of fluid to bounce on a free surface has been
observed in different contexts. \cite{BounceJet} designed an
experiment where a jet of Newtonian fluid falls into a rotating vessel
filled with the same fluid. The authors investigated the conditions
under which the jet bounces: nature of the fluid, jet diameter, jet and
bath velocities.  We refer to \cite{lockhart} for experimental movies
illustrating this effect.
\begin{figure}[h]
\centering
\includegraphics[scale=0.2]{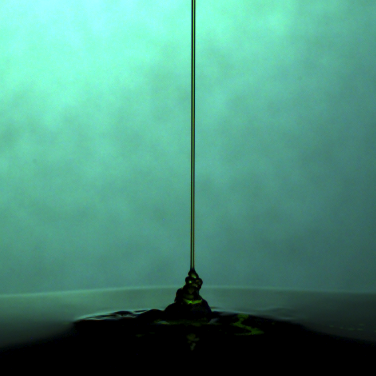}
\includegraphics[scale=0.2]{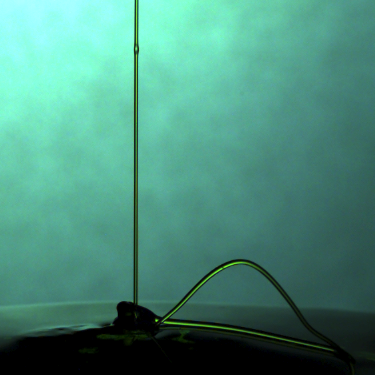}
\includegraphics[scale=0.2]{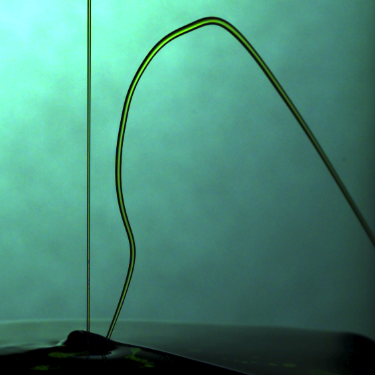}
\caption{Experimental Observation of the Kaye Effect. (Left) the fluid start to buckle producing a heap; (middle) A stream of liquid suddenly leaps outside the heap; (right) Fully developed Kaye effect.}
\label{fig:Kaye_effect}
\end{figure}
Bouncing can also be observed in a stationary vessel provided the jet
is composed of a \emph{non-newtonian} fluid.  During the pouring
process, a small heap of fluid forms and the jet occasionally leaps
upward from the heap; see Figure \ref{fig:Kaye_effect}.  This is the
so-called \emph{Kaye} effect as first observed by \cite{K_1963} in
1963.  About 13 years after this phenomenon was first mentioned in the
literature, \cite{C_1976} revisited the experiment and suggested that
the ability for the fluid to exhibit shear-thinning viscosity and
elastic behavior are key ingredients for the bouncing to occur.
Additional laboratory experiments performed by
\cite{versluis2006leaping} and \cite{B_2009} lead the authors of each
team to propose a list of properties that the fluid should have for
the Kaye effect to occur.  The conclusions of these two papers
disagree on the requirement that the fluid be elastic and on the
nature of the thin layer that separates the heap and the outgoing jet.
It is argued in \cite{versluis2006leaping} that elastic properties are
not necessary and that the thin layer is a shear layer, whereas it is
argued in \cite{B_2009} that the fluid should have elastic properties
and that the thin layer separating the heap from the bouncing jet is a
layer of air.  Recent experiments reported in
\cite{PhysRevE.87.061001} using a high speed camera unambiguously show
that the jet slides on a lubricating layer of air.  For completeness,
we also refer to \cite{Klapp_2012} for a thorough discussion on the
``stable'' Kaye effect, where the jet falls against an inclined
surface.

The objective of the present paper is to numerically revisit the Kaye
effect.  Our key finding is that \emph{a thin air layer is always
  present between the bouncing jet and the rest of the fluid whether
  the fluid is Newtonian or not.} Our numerical experiments suggest
that the critical parameter for bouncing to occur is that the
properties of fluid and the experimental conditions be such that a
stable layer of air separating the jet and the ambiant fluid can
appear.  This condition is met by setting the bath in motion for
Newtonian fluids; it can also be met if the bath is stationary
provided the fluid is non-Newtonian and has shear-thinning
viscosity. The numerical code that we use is based on a modeling 
the fluid flows by the incompressible Navier-Stokes equations
supplemented with a surface tension mechanism.  The shear-thinning
viscosity of the non-Newtonian fluids is assumed to follow a model by
\cite{Cross1965417}.  Elastic behaviors are not modeled. The numerical
approximation of the resulting two-phase flow model is based on two
solvers: one solving the Navier-Stokes equations assuming that the
fluid/air distribution is given, the other keeping track of the motion
of the interface assuming the transport velocity is given. The
Navier-Stokes solver is based on projection method by
\cite{Chorin1968} and \cite{Temam1969} using a second-order backward
differentiation formula for the time discretization and finite
elements for the space approximation.  The transport solver is based
on a level-set technique in the spirit of \cite{osher1988fronts} to
represent the liquid/air interface. The level set is approximated in
space by using finite elements and the time stepping is done by using
a third order explicit Runge-Kutta technique.

The paper is organized as follows. The mathematical model is presented
in Section \ref{sec:math_model}.  The numerical techniques to solve
the Navier-Stokes equations and the transport equation for the level
set function are described in Section \ref{sec:num_method}.  Various
validation tests of the numerical algorithms and comparisons with
classical benchmark problems are reported in Section
\ref{s:validation}.  Finally, we report numerical evidences of
Newtonian and non-Newtonian bouncing jets in Section
\ref{sec:bouncing_jets}.


\section{The Mathematical Model}
\label{sec:math_model}
This section presents the mathematical models adopted to describe
non-mixing two phase fluid flows with capillary forces.  Each fluid is
assumed to be incompressible.

\subsection{Two Phase Flow System}
\label{subsec:two_phase}
Let $\Lambda \subset \mathbb{R}^d$ ($d=2,3$) be an open and bounded
computational domain with Lipschitz boundary $\partial \Lambda$ and
let $[0,T]$ be the computational time interval, $T>0$.  The cavity
$\Lambda$ is filled with two non-mixing fluids undergoing some
time-dependent motion, say fluid 1 and fluid 2. We denote by
$\Omega^+$ and $\Omega^-$ the open subsets of the space-time domain
$\Lambda\CROSS [0,T]$ occupied by fluid 1 and fluid 2, respectively.
We denote by $\rho^+$, $\mu^+$, $\rho^-$, $\mu^-$ the density and
dynamical viscosity of each fluid, respectively. The interface between
the two fluids in the space-time domain is denoted $\Sigma :=
\partial\Omega^+\cap \partial\Omega^-$, and the normal to $\Sigma$, oriented from
$\Omega^+$ to $\Omega^-$, is denoted $\bn_\Sigma$. The two space-time
components of the vector field $\bn_\Sigma : \Sigma \longrightarrow
\Real^d \CROSS \Real$ are denoted $\bn$ and $\bn_\tau$,
respectively. It is also useful to define $\Omega^\pm(t) :=
\Omega^\pm\cap (\Lambda\CROSS\{t\})$; \ie the sets $\Omega^+(t)$ and
$\Omega^-(t)$ are the regions occupied by fluid 1 and fluid 2 at time
$t$, respectively. We also introduce the interface $\Sigma(t) =
\partial\Omega^+(t)\cap \partial\Omega^-(t)$; note that $\bn$, as defined above, is
the unit normal of $\Sigma(t)$ and it is oriented from $\Omega^+(t)$ to
$\Omega^-(t)$.  To facilitate the modeling, we define global density
and dynamical viscosity functions $\rho, \mu : \Lambda \CROSS [0,T]
\rightarrow \mathbb R$ by setting $\rho(\bx,t) = \rho^\pm$ and
$\mu(\bx,t) = \mu^\pm$ if $(\bx,t)\in \Omega^\pm$. The fluid velocity
field ${\bu} :\Lambda \CROSS [0,T] \rightarrow \mathbb{R}^d$,
henceforth assumed to be continuous across $\Sigma$, and the pressure
$p: \Lambda \CROSS [0,T] \rightarrow \mathbb{R}$ are defined globally
and solve the incompressible Navier-Stokes equations in the
distribution sense in the space-time domain
\begin{subequations}\label{eqn:navier_stokes}
\begin{align}
  \rho \left( \dfrac{\partial}{\partial t} \bu + \bu \SCAL \nabla
    \bu \right) - 2 \mbox{div}\left( \mu\ \nabla^S\bu \right) + \nabla
  p - \delta_{\Sigma}\sigma \kappa \bn
  &= \rho \bg   & \mbox{ in } &\Lambda \CROSS (0,T],  \label{eqn:navier_stokes_mom}  \\
  \mathrm{div} (\bu)&=0 & \mbox{ in } &\Lambda \CROSS
  (0,T], \label{eqn:navier_stokes_mass}
\end{align}
\end{subequations}
where $\nabla^S\bu := \frac 1 2 (\nabla \bu + \nabla \bu^T)$ is the
strain rate tensor, $\bg$ is the gravity field, and
$\delta_{\Sigma}\sigma \kappa \bn$ is a singular
measure modeling the surface tension acting on $\Sigma(t)$. The
distribution $\delta_{\Sigma}$ is the Dirac measure supported on
$\Sigma$, the function $\kappa:\Sigma \rightarrow \mathbb R$ is
the total curvature of $\Sigma(t)$ (sum of the principal curvatures)
and $\sigma$ is the surface
tension coefficient.  

The system \eqref{eqn:navier_stokes} is supplemented with initial and
boundary conditions.  The initial condition is $\bu(\bx,0)=\bu_0(\bx)$
for all $\bx\in \Lambda$, where $\bu_0$ is assumed to be a smooth
divergence-free velocity field.  The boundary $\partial \Lambda$ is
decomposed into three non-overlapping components $\partial
\Lambda:=\overline{\Gamma}_D\cup \overline{\Gamma}_N
\cup\overline{\Gamma}_N$ with $\Gamma_D\cap\Gamma_N = \emptyset$,
$\Gamma_N\cap\Gamma_S = \emptyset$, $\Gamma_S\cap\Gamma_D =
\emptyset$.  Time-dependent decomposition could be considered but are
not described here to avoid unnecessary technicalities. Given
$\bef_N:\Gamma_N\rightarrow \mathbb R^d$ and
$\bef_D:\Gamma_D\rightarrow \mathbb R^d$, we require that
\begin{subequations}
\label{e:navier_stokes_bdy_cond}
\begin{align}
\left(2\mu \nabla^S\bu - pI\right)\bnu &= \bef_N, & \text{ on } &\Gamma_N\CROSS(0,T],\\
\bu &= \bef_D, & \text{ on } &\Gamma_D \CROSS (0,T], \\
\bu \SCAL \bnu =0, \quad \left( \left(2\mu \nabla^S\bu - pI\right)\bnu\right)\CROSS \bnu&=0, &  
\text{ on } &\Gamma_S \CROSS (0,T], \label{eqn:slip}
\end{align}
\end{subequations}
where $\bnu$ is the outward pointing unit normal on $\partial\Lambda$ and $I$
is the $d\CROSS d$ identity matrix.  For simplicity, we assume that
the $(d-1)$-measures of $\Gamma_N\cup \Gamma_S$ and $\Gamma_D$ are
each strictly positive; otherwise extra constraints either on the
velocity or on the pressure must be enforced.

Note that we could have formulated the
conservation of momentum without invoking the singular measure
modeling the surface tension by saying that \eqref{eqn:navier_stokes_mom}
holds in $\Omega^+$ and $\Omega^-$ (without the singular measure) and
by additionally requiring that
\begin{equation}\label{e:interface} [ \bu ] = 0 \quad \text{and} \quad
  \left[2\mu \nabla^S\bu - pI\right]\bn = \sigma \kappa \bn \quad
  \text{on} \quad \Sigma(t),\qquad \forall t\in (0,T],
\end{equation}
where $[.]$ denotes the jump across $\Sigma(t)$ defined by
$[\bv](\bx):= \lim_{\Omega^-\ni\by\to \bx} \bv(\by) -
\lim_{\Omega^+\ni\bz\to \bx} \bv(\bz) $, \ie $[\bv](\bx) = \bv(\bx^-)
- \bv(\bx^+)$, 
for all $\bx \in
\Sigma(t)$ and all $\bv:\Lambda \rightarrow \mathbb R^d$ or
$\bv:\Lambda \rightarrow \mathbb R^{d\CROSS d}$.

The interface $\Sigma(t)$ is assumed to be transported by the fluid
particles. More precisely, let $\partial^-(\Lambda\CROSS [0,T])$ be the
part of the boundary of the space-time domain where the
characteristics generated by the field $(\bu,1)$ enter, \ie
\begin{equation}
  \partial^-(\Lambda\CROSS [0,T]) := \Lambda
  \CROSS\{0\} \cup \{(\bx,t) \in \partial\Lambda\CROSS[0,T]\ | \
  \bu(\bx,t)\SCAL\bnu < 0\}.
\end{equation}  
We then define
$\{\bx(\bP,t)\in \Lambda, \ t\in [t_\bP,T],\
(\bP,t_\bP)\in \partial^-(\Lambda\CROSS [0,T])\}$
to be the family of the characteristics generated by the velocity
field $\bu$, \ie
$\frac{\partial}{\partial t}\bx(\bP,t) = \bu(\bx(\bP,t),t)$ with
$\bx(\bP,t_\bP)=\bP$, $(\bP,t_\bP)\in\partial^-(\Lambda\CROSS [0,T])$
where $t_\bP$ is the time when the characteristics enters the
space-time domain $\Lambda\CROSS [0,T]$ at $\bP$.  Let us now denote
by $\Sigma_0:= \Sigma\cap \partial^-(\Lambda\CROSS [0,T])$ the
location of the interface at the inflow boundary of the space-time
domain, then we are going to assume in the entire paper that the
velocity field is smooth enough so that the following property holds
\begin{equation}\label{eqn:free_surface}
\forall \bx\in \Sigma(t), \ \exists! (\bP,t_\bP)\in \Sigma_0, \ 
\bx = \bx(\bP,t),\ t\ge t_\bP.
\end{equation}

\subsection{Eulerian Representation of the Free Boundary Interface}
A level set technique is used to keep track of the position of the
time-dependent interface $\Sigma(t)$, see for instance
\cite{osher1988fronts}.  This method is recalled in Section
\ref{subsec:level_set}. A typical problem arising when using a
level-set method to describe interfaces is to guarantee the
non-degeneracy of the representation.  We discuss a reinitialization
technique overcoming this issue in Section
\ref{subsec:levelset:reinit}.

\subsubsection{Level-Set Representation}\label{subsec:level_set}
Let us define the so-called level function $\phi(\bx,t) : \Lambda
\CROSS [0,T] \rightarrow \mathbb{R}$ so that 
\begin{equation}
  \frac{\partial}{\partial t}\phi + \bu\SCAL\nabla \phi = 0 
  \quad \text{in} \quad \Lambda \CROSS (0,T],
  \qquad \phi(\bP,t_\bP) =\phi_0(\bP,t_\bP)\quad \text{on} 
  \quad \partial^-(\Lambda\CROSS [0,T]),
\label{eqn:level_set}
\end{equation}
where we assume that $\phi_0$ is a smooth function satisfying the
following properties:
\begin{align}
\partial\Omega^\pm \cap\partial^-(\Lambda\CROSS[0,T])
&= \{ (\bP,t_{\bP}) \in\partial^-(\Lambda\CROSS[0,T])  \ | \ \pm\phi_0(\bP,t_\bP) \ge 0 \}.\\
\Sigma_0:=\Sigma \cap\partial^-(\Lambda\CROSS[0,T])
&= \{ (\bP,t_{\bP}) \in\partial^-(\Lambda\CROSS[0,T])  \ |\ \phi_0(\bP,t_\bP) = 0 \}.
\end{align}
Note that this definition implies that $\Sigma_0$ is the $0$-level set
of $\phi_0$.  Upon introducing
$\psi(\bP,t):=\phi(\bx(\bP,t),t)$ for all
$(\bP,t_\bP)\in \partial^-(\Lambda\CROSS[0,T])$, $t\ge t_\bP$, the
definition of $\phi$ together with the definition of the
characteristics implies that $\partial_t\psi(\bP,t) = 0$; thereby
proving that $\phi(\bx(\bP,t),t) = \psi(\bP,t) = \psi(\bP,t_\bp) =
\phi(\bx(\bP,t_\bP),t_\bP) = \phi(\bP,t_\bP) =\phi_0(\bP,t_\bP)$. This
means that the value of $\phi(\bx(\bP,t),t)$ along the trajectory
$\{\bx(\bP,t_\bP),\ t\in [t_\bP,T]\}$ is constant; in particular the
sign of $\phi(\bx(\bP,t),t)$ does not change.  The leads us to adopt
the following alternative definitions for $\Omega^\pm$ and
$\Sigma(t)$:
\begin{subequations}
\label{Levelset_Characterization}
\begin{align}
\Omega^\pm &:=\{(\bx,t) \in \Lambda\CROSS [0,T]\ | \pm\phi(\bx,t)>0\} \\
\Sigma(t) &:=  \{\bx \in \Lambda\ | \ \phi(\bx,t)=0\}.
\end{align}
\end{subequations}
The above characterizations will be used in the rest of the paper.

To simplify the presentation we are going to assume in the rest of the
paper that there is $\Gamma_L \subset \partial\Lambda$ such that
$\partial^-(\Lambda\CROSS(0,T))= \Gamma_L\CROSS (0,T)$, \ie the inflow
boundary for the level-set equation is time-independent.  We then set
$\phi_L(\bP,t) = \phi_0(\bP,t)$ for all $\bP\in \Gamma_L$, $t\in
[0,T]$, and we abuse the notation by setting $\phi_0(\bP) =
\phi_0(\bP,0)$ for all $\bp\in \Lambda$.

\subsubsection{Reinitialization and cut-off function}
\label{subsec:levelset:reinit}
A typical issue when dealing with level-set representation of
interfaces is to guarantee that the manifold $\{\bx \in \Lambda\ | \
\phi(\bx,t)=0\}$ is $(d-1)$-dimensional, \ie we want to make sure that
$\|\nabla \phi \|_{\ell^2}>0$ in every small neighborhood of the zero
level-set, where $\|\cdot\|_{\ell^2}$ is the Euclidean norm in
$\Real^d$. In order to achieve this objective, we implement an ``on
the fly'' reinitialization algorithm proposed in
\cite{ville2011convected}, which consists of replacing
\eqref{eqn:level_set} by
\begin{equation}
  \dfrac{\partial}{\partial{t}}{\phi_{\lambda,\beta}} 
+\bu \SCAL \nabla \phi_{\lambda,\beta}
= \lambda \mbox{sign}(\phi_{\lambda,\beta}) \left( \mathcal G(\phi_{\lambda,\beta})- \|\nabla \phi_{\lambda,\beta}\|_{\ell^2} \right),
\label{eqn:level:final}
\end{equation}
where $\lambda,\beta>0$ are parameters yet to be defined, and
$\mathcal G(\cdot)$, $\mbox{sign}(\cdot)$ are defined by
\begin{equation}\label{eqn:level:sign}
\mathcal G(z) :=1 - \left(\frac{z}{\beta} \right)^2 ,\qquad 
\text{sign}(z) :=
\left\{ 
{\begin{array}{rl}
1,  & z   > 0 \\
-1, & z   < 0 \\ 
0,  & |z| = 0.
\end{array}} 
\right.
\end{equation} 
The rational for the new definition \eqref{eqn:level:final} is that
the presence of the sign function in the right-hand side implies that
the $0$-level set of $\phi_{\lambda,\beta}$ is the same as that of
$\phi$; \ie the characterizations of $\Omega^\pm$ and $\Sigma(t)$ are
unchanged, see \eqref{Levelset_Characterization}.  Moreover, assuming
that $\bu$ is locally the velocity of a solid motion and upon setting
$\psi(\bP,t) = \phi_{\lambda,\beta}(\bx(\bP,t))$, we have
$\frac{\partial}{\partial{t}}{\psi} = \lambda \mbox{sign}(\psi) ( (1 -
(\frac{\psi}{\beta})^2) - |\nabla \psi|)$. Assuming that this eikonal
equation has a steady state solution $\psi_\infty$ and denoting by
$\Sigma_\infty$ the $0$-level set of this steady state solution, the
behavior of $\psi_\infty$ in the vicinity of $\Sigma_\infty$ is
$\psi_\infty(\bp) \approx \beta
\tanh\left(\frac{\text{dist}(\Sigma_\infty,\bP)}{\beta}\right)$, since
the solution to the following $ODE$: $y'(z) = 1 - (y(z)/\beta)^2$ is
\begin{equation}
\mathcal F(z):= \beta \tanh\left( \dfrac{z}{\beta}\right).
\label{eqn:cut_off_beta}
\end{equation}
Note that $\psi$ is close to $\psi_\infty$ if $\lambda$ is very large.
In conclusion the solution of \eqref{eqn:level:final} is such that
$\phi_{\lambda,\beta}(\bx,t) \approx \beta
\tanh\left(\frac{\text{dist}(\Sigma(t),\bx)}{\beta}\right)$ if
$\lambda$ is large enough, \ie $\|\nabla \phi\|_{\ell^2} \sim 1$ in any small
neighborhood of the zero level-set. We are going to abuse the notation
in the rest of the paper by dropping the indices $\lambda$, $\beta$
and by using $\phi$ instead of $\phi_{\lambda,\beta}$.

\section{Numerical Method}
\label{sec:num_method}
We discuss the approximations of the level-set equation in
Section~\ref{subsec:num_level_set} and that of the Navier-Stokes
system in Section~\ref{ss:NS}.  We consider a mesh family,
$\{\mathcal T_h\}_{h>0}$, and we assume that each mesh $\mathcal T_h$
is a subdivision of $\overline{\Lambda}$ made of disjoint elements
$K$, \ie rectangles when $d=2$ or cuboids when $d=3$. We denote by
$\mathcal E(\mathcal T_h)$ the collection of interfaces and
boundary faces (edge when $d=2$ and faces when $d=3$).  Each
subdivision is assumed to exactly approximate the computational
domain, \ie $\overline{\Lambda} = \cup_{K \in \mathcal{T}_h} K$, and
to be consistent with the decomposition of the boundary, \ie there
exists $\mathcal E_I(\mathcal T_h) \subset \mathcal E(\mathcal T_h)$
for $I \in \{D,N,S,L\}$ such that
$\overline{\Gamma}_I = \cup_{F \in \mathcal E_I(\mathcal T)} {F}$.
The diameter of an element $K \in \mathcal T_h$ is denoted by $h_K$;
$h_K:=\max_{\bx,\by\in K} |\bx-\by|$.  The mesh family
$\{\mathcal T_h\}_{h>0}$ is assumed to be shape regular in the sense
of Ciarlet.  For any integer $k \geq 1$ and any $K \in \mathcal{T}_h$,
we denote by $\mathbb{Q}^k(K)$ the space of scalar-valued multivariate
polynomials over $K$ of partial degree at most $k$.  The vector-valued
version of $\mathbb{Q}^k(K)$ is denoted $\pmb{\mathbb Q}^k(K)$.  The
index $h$ is dropped in the rest of the paper and we write
$\mathcal T$ instead of $\mathcal T_h$ when the context is
unambiguous.

Regarding the time discretization, given an integer $N\geq 2$, we
define a partition of the time interval $0=:t^0 < t^1 < ... < t^N:=T$
and denote $\delta t^n:=t^n-t^{n-1}$ and $t^{n+\frac 1 2}:=
t^{n}+\frac 1 2 \delta t^{n+1}$.

\subsection{Numerical Approximation of the Level-Set System}
\label{subsec:num_level_set}
The continuous finite element method used for the space approximation
of the level-set equation \eqref{eqn:level:final} is described in
Section~\ref{eqn:levelset:apprx_space}.  We present in Sections
\ref{sec:levelset_entropy} an entropy-residual technique that has the
advantage of avoiding the spurious oscillations that would otherwise
be generated by using an un-stabilized Galerkin technique.  The time
stepping is done by using an explicit strong stability preserving
(SSP) Runge-Kutta 3 (RK3) scheme as explained in Section
\ref{eqn:level set:time}.

\subsubsection{Approximation in Space}\label{eqn:levelset:apprx_space}
The space approximation of $\Phi(\cdot,t)$ of the level set function $\phi(\cdot,t)$ solution
to \eqref{eqn:level:final} is done by using continuous, piecewise linear
polynomials subordinate to the subdivision $\mathcal T$. 
The associated finite element spaces are
defined by
\begin{align}
\mathbb{W}(\mathcal{T}) &:= \{ W \in C^0(\overline{\Lambda};\mathbb{R}) \ | 
\ W |_{K} \in \mathbb{Q}^1(K) , \ \forall K \in \mathcal{T} \},\\
\mathbb{W}_0(\mathcal{T}) &:= \{ W \in C^0(\overline{\Lambda};\mathbb{R})\ | 
 \ W= 0 \text{ on } \Gamma_{L}, \ W |_{K} \in \mathbb{Q}^1(K) , \ \forall K \in \mathcal{T} \},\\
\mathbb{W}_{L}(\mathcal{T}) &:= \{ W \in C^0(\overline{\Lambda};\mathbb{R}) \ | 
\ W= \Phi_L \text{ on } \Gamma_L,
\ W |_{K} \in \mathbb{Q}^1(K) , \ \forall K \in \mathcal{T} \},
\end{align}
where $\Phi_L(\cdot,t)$ is a piecewise linear approximation of the
inflow data $\phi_L(\cdot,t)$.  Assuming that the velocity field $\bu:
\Lambda\CROSS[0,T]\longrightarrow \Real^d$ is known, the Galerkin
approximation of \eqref{eqn:level:final} is formulated as follows:
Given $\Phi_L$ and $\Phi(\cdot,0) = \Phi_0$, where $\Phi_0 \in \mathbb{W}(\mathcal{T})$
is an approximation of the initial condition $\phi_0$, find $\Phi \in \mathcal
C^1([0,T];\mathbb{W}_L(\mathcal{T}))$ such that
\begin{equation}
\int_\Lambda \dfrac{\partial}{\partial{t}}\Phi  W
= - \int_\Lambda  \bu \SCAL \nabla \Phi \, W  
+ \int_\Lambda \lambda \text{sign}(\Phi) \left(\mathcal G(\Phi) - \|\nabla\Phi\|_{\ell^2}\right)  W
 , \quad \forall W \in \mathbb{W}_0(\mathcal{T}).
\label{eqn:level_set:space_h_Gal}
\end{equation} 
It is well known that the solution to the above system exhibit
spurious oscillations in the regions where $\|\nabla \Phi\|_{\ell^2}$ is large.
We address this issue in the next section.

\subsubsection{Entropy Residual Stabilization}

\label{sec:levelset_entropy}
\begin{figure}[h]
\centering
\begin{subfigure}[]{0.35\textwidth}
\captionsetup{justification=centering}
\includegraphics[width=\textwidth]{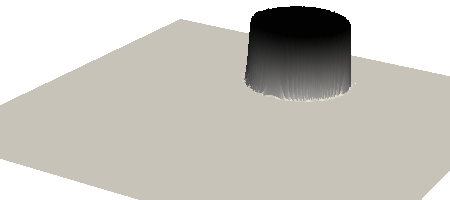}
\caption{}\label{fig:level:fv_0_a}
\end{subfigure}
\begin{subfigure}[]{0.35\textwidth}
\captionsetup{justification=centering}
\includegraphics[width=\textwidth]{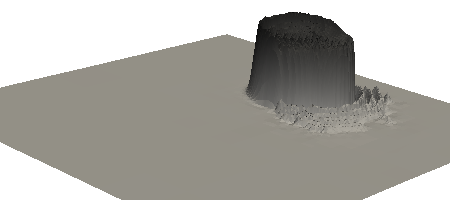} 
\caption{}
\label{fig:level:fv_0_b}
\end{subfigure}
\begin{subfigure}[]{0.35\textwidth}
\captionsetup{justification=centering}
\includegraphics[width=\textwidth]{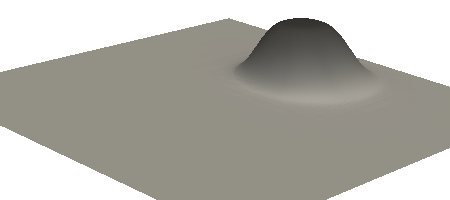} 
\caption{}
\label{fig:level:fv_0_c}
\end{subfigure}
\begin{subfigure}[]{0.35\textwidth}
\captionsetup{justification=centering}
\includegraphics[width=\textwidth]{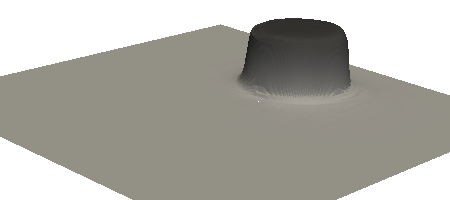} 
\caption{}
 \label{fig:level:fv_0_d}
\end{subfigure}
\caption{Graph of the level-set function in the computational domain
  $\Lambda = (-1,1)^2$ using the velocity given by
e  \eqref{eqn:level:circle_vel}: (a) Initial data: $\phi_0$ is the
  characteristic function of the disk of radius $1$ centered at
  $(\frac12,0)$; (b) No stabilization; (c) First-order stabilization
  with $C_{\text{Lin}} = 0.1$; (d) Entropy viscosity stabilization
  with $C_{\text{Lin}} = 0.1$ and $C_\text{Ent} = 0.1$. Observe the
  spurious oscillations in panel (b) when no stabilization is applied.
  Both, the first-order and the entropy viscosity solution are free of
  oscillations, the latter being clearly more accurate.}
\label{fig:level:fv_0}
\end{figure}

We describe in this section an entropy viscosity technique to
stabilize the Galerkin formulation \eqref{eqn:level_set:space_h_Gal}. This
method has been introduced in \cite{guermond2011entropy} and we refer
to \cite{MR3167449} for a mathematical discussion on its stability
properties.  To motivate the discussion, we refer to the panel (b) of
Figure \ref{fig:level:fv_0} showing the Galerkin approximation of the
characteristic function of the unit disk, initially centered at
$(\frac12,0)$, after one rotation about the origin.

The spurious oscillations are avoided by augmenting
\eqref{eqn:level_set:space_h_Gal} with an artificial viscosity term where
the viscosity is localized and chosen to be proportional to an entropy
residual. To describe the method and define an appropriate local
``viscosity", we recall that the following holds in the distribution sense 
for any $E\in \calC^1(\Lambda;\Real)$:
\[
\dfrac{\partial}{\partial t} E(\phi) +   \bu\SCAL \nabla E(\phi)  
-  \lambda \mbox{sign}(\phi) \left(\mathcal G(\Phi)- \|\nabla \phi\|_{\ell^2} \right)E'(\phi) = 0,
\]
Consequently, it is reasonable to expect that the semi-discrete entropy residual
\begin{equation*}
  R^{\mathrm{Ent}}(\Phi,\bu) := \dfrac{\partial}{\partial t} E(\Phi) 
  +   \bu(t) \SCAL \nabla  E(\Phi) 
-  \lambda \mbox{sign}(\Phi) \left(\mathcal G(\Phi)- \|\nabla \Phi\|_{\ell^2} \right)E'(\Phi)
  b\end{equation*}
is a reliable indicator of the regularity of $\phi$.  This quantity
should be of the order of the consistency error in the regions where
$\phi$ is smooth and it should be large in the region where the PDE is
not well solved.  In our computations, we have chosen 
\begin{equation}\label{eqn:entropy:choice}
E(\phi) = |\phi|^p,
\end{equation}
The local so-called entropy viscosity is defined for any $K\in \mathcal T$ by
\begin{equation}
{\mu}^{\text{Ent}}_{K}(\Phi,\bu) := C_{\text{Ent}}h_K^2  \frac{\| R^{\text{Ent}}(\Phi,\bu) \|_{L^\infty(K)}}
{\| E(\Phi)- \frac 1 {|\Lambda|} \int_\Lambda E(\Phi) \|_{L^\infty(\Lambda)}},
\label{eqn:level:fv}
\end{equation}
where $C_{\text{Ent}}$ is an absolute constant.  In the regions where
$\phi$ is discontinuous (or has a very sharp gradient), the entropy
viscosity as defined above may be too large and thereby introduce too
much diffusion, which in turn may severely limit the CFL number when
using an explicit time stepping.  In this case a linear first-order
viscosity is turned on instead
\begin{equation}\label{eqn:lin:visc}
 \mu^{\text{Lin}}_{K}(\Phi,\bu) = C_{\text{Lin}} 
\|h_K ({\bu}+ \lambda \mbox{sign}(\Phi) \frac{\nabla \Phi }{\|\nabla \Phi\|_{\ell^2}} ) \|_{L^{\infty}(K)},
\end{equation}
where $C_{\text{Lin}}$ is an absolute constant. The justification for
the definition of the local speed that is used to define the viscosity
in \eqref{eqn:lin:visc} is that \eqref{eqn:level:final} can be
re-written $\frac{\partial}{\partial{t}}{\phi} +\bw \SCAL \nabla \phi
= \lambda \mbox{sign}(\phi) \mathcal G(\phi)$, where $\bw = \bu +
\lambda \mbox{sign}(\phi) \frac{\nabla \phi }{\|\nabla \phi\|_{\ell^2}}$.
Combining the two viscosities yield the artificial viscosity
$\mu^{\text{Stab}}:\Lambda \CROSS [0,T] \rightarrow \mathbb R$ defined
on each $K\in \mathcal T$ by
\begin{equation}\label{eqn:levelset_stab_v}
{\mu}^{\text{Stab}}(\Phi,\bu)|_K := \min({\mu}^{\text{Lin}}_{K}(\Phi,\bu) , \mu^{\text{Ent}}_{K}(\Phi,\bu) ).
\end{equation}

Going back to the space discretization, we modify
\eqref{eqn:level_set:space_h_Gal} as follows: Look for $\Phi\in
\mathcal C^1([0,T];\mathbb W_L(\mathcal T))$ so that
\begin{equation}
\int_\Lambda \dfrac{\partial}{\partial{t}}\Phi W
= - \int_\Lambda  \bu \SCAL \nabla \Phi \, W  
+ \int_\Lambda \lambda \text{sign}(\Phi) \left(\mathcal G(\Phi)- \|\nabla\Phi\|_{\ell^2}\right)  W
- \int_\Lambda {\mu}^{\text{Stab}}(\Phi,\bu)\nabla \Phi \SCAL \nabla W 
\label{eqn:level_set:space_h}
\end{equation}
for all $W \in \mathbb{W}_0(\mathcal{T})$
and $\Phi(\cdot,0) = \Phi_0$.

\subsubsection{Approximation in Time}\label{eqn:level set:time}
Before introducing the time discretization we re-write
\eqref{eqn:level_set:space_h} as follows:
\begin{equation}
\int_\Lambda \dfrac{\partial}{\partial t} \Phi \, W = \int_\Lambda L(\Phi,\bu)W - \int_\Lambda \mu^{\textrm{Stab}}(\Phi,\bu) \nabla \Phi \SCAL \nabla W, 
\qquad  \forall W \in \mathbb{W}_0(\mathcal{T},t), 
\label{eqn:level:time}
\end{equation}
where $L(\Phi,\bu) := - \bu + \lambda \text{sign}(\Phi) \left(\mathcal
  G(\Phi) -\|\nabla \Phi\|_{\ell^2}\right)$. Then we approximate time in the
above nonlinear system of ODES by using an explicit RK3 strong
stability preserving (SSP) scheme, \eg see
\cite{gottlieb2001strong,Shu_RK3} for more details on SSP methods.  We
denote by $\Phi^k$ the approximation of $\Phi(\cdot,t^k)$, $0\le k\le
N$.  Then, the time stepping proceeds as follows: Given $\Phi^{n}$,
compute $\Phi^1$, $\Phi^2, \Phi^3 \in W(\mathcal T)$ and $\Phi^{n+1}
\in \mathbb W_L(\mathcal T)$ so that
\begin{align} \label{eqn:level_ssprk3}
(i) \qquad &\int_\Lambda \Phi^{(1)} \, W = \int_\Lambda \left( \Phi^n\ + \delta t^{n+1} L(\Phi^n,\bu^n) \right) W, \qquad \forall W \in \mathbb{W}(\mathcal{T})  \nonumber \\ 
(ii) \qquad&\int_\Lambda \Phi^{(2)} \, W = \int_\Lambda \left( \dfrac{3}{4}\Phi^n + \dfrac{1}{4} \Phi^{(1)} +\dfrac{1}{4} \delta t^{n+1} L(\Phi^{(1)},\bu^{n+1}) \right) W - \int_\Lambda \mu^{\mathrm{Stab},n+1} \nabla \Phi^{(1)} \SCAL \nabla W, \qquad \forall W \in \mathbb{W}(\mathcal{T})   \\ 
(iii) \qquad&\int_\Lambda \Phi^{(3)} \, W = \int_\Lambda \left( \dfrac{1}{3}\Phi^n + \dfrac{2}{3} \Phi^{(2)} +\dfrac{2}{3} \delta t^{n+1} L(\Phi^{(2)},\bu^{n+\frac12}) \right) W- \int_\Lambda \mu^{\mathrm{Stab},n+\frac12} \nabla \Phi^{(2)} \SCAL \nabla W, \qquad \forall W \in \mathbb{W}(\mathcal{T})  \nonumber
\end{align}
and 
\begin{equation}\label{eqn:level_ssprk3_2}
\Phi^{n+1}(\bx,t) = \left\lbrace \begin{array}{ll}
 \Phi^{(3)}(\bx) & \quad \textrm{when  } \bx \not \in \Gamma_L,\\
 \Phi_L(\bx,t^{n+1}) &  \quad \textrm{when  } \bx  \in \Gamma_L.
 \end{array}\right.
 \end{equation}
 Using \eqref{eqn:levelset_stab_v}, the viscosities
 $\mu^{\mathrm{Stab},n+1}$ and $\mu^{\mathrm{Stab},n+\frac12}$ are
 defined to be equal to $\mu^{\textrm{Stab}}(\Phi^{(1)},\bu^{n+1})$
 and $\mu^{\textrm{Stab}}(\Phi^{(2)},\bu^{n+\frac12})$, respectively,
 where the residual $R^{\textrm{Ent}}$ in the definition
 \eqref{eqn:level:fv} of $\mu^{\textrm{Ent}}$ is evaluated as follows:
\begin{align}\label{eqn:Residual1}
R(\Phi^{(1)},\bu^{n+1})  &\approx \frac{E(\Phi^{(1)})-E(\Phi^n)}{\delta t^{n+1}} 
+ \left(\bu(t^{n+1}) \SCAL \nabla \Phi^{(1)} 
-\lambda \mbox{sign}(\Phi^{(1)}) 
\left(\mathcal G(\Phi^{(1)})-\|\nabla \Phi^{(1)}\|_{\ell^2}\right)\right)E'(\Phi^{(1)})\\
\label{eqn:Residual2}
R(\Phi^{(2)},\bu^{n+\frac 12}) &\approx
\frac{E(\Phi^{(2)})-E(\Phi^n)}{\frac 1 2 \delta t^{n+1}} +
\left(\bu(t^{n+\frac12}) \SCAL \nabla \Phi^{(2)} -\lambda
  \mbox{sign}(\Phi^{(2)}) \left(\mathcal G(\Phi^{(2)})-\|\nabla
    \Phi^{(2)}\|_{\ell^2}\right)\right)E'(\Phi^{(2)}).
\end{align}
\begin{rmk}[``On the Fly" Stabilization]
  Notice that following \cite{MR3167449}, no viscosity is added in the
  computation of $\Phi^{(1)}$.  In particular, the viscosities used
  within the time interval $[t^n,t^{n+1}]$ only depend on the values
  of $\Phi$ and $\bu$ on the same time interval.
\end{rmk}
\begin{rmk}[Stability and Convergence]
  We expect the scheme to be stable under the following
  Courant-Friedrichs-Lewy (CFL) condition:
\begin{equation}
\delta t^{n+1} \leq C_{CFL}\min_{K \in \mathcal{T}} \frac{h_K}{\|\bu(t^{n+1})
+  \lambda \mbox{\rm sign}(\Phi^n) \frac{\nabla \Phi^n}{\|\nabla \Phi^n\|_{\ell^2}}\ \|_{L^{\infty}(K)}},
\label{eqn:CFL}
\end{equation}
for some sufficiently small but positive constant $C_{CFL}$
independent of $\lambda$, $\mathcal{T}$, $\delta t$, $\Phi$ and $\bu$.
Refer for instance to \cite{MR3167449} for further details on the CFL
condition.
Moreover, it seems reasonable to expect that only the entropy
viscosity is active in the regions where $\phi$ is smooth; as a
result, the CFL condition implies that the first-order approximation
of the time derivative in the evaluation of the entropy residuals,
\eqref{eqn:Residual1}-\eqref{eqn:Residual2}, does not affect the
overall third-order approximation thanks to the $h_K^2$ factor present
in the definition of the viscosity \eqref{eqn:level:fv}.  This
conjecture is confirmed numerically in Section
\ref{sec:validation_transport}.
\end{rmk}


\subsection{Numerical Approximation of The Navier-Stokes System}
\label{ss:NS}
The space approximation of the velocity and pressure in the
Navier-Stokes equations is done by using Taylor-Hood finite
elements. The time discretization is done by using the second-order
backward differentiation formula (BDF2).  An incremental rotational
pressure correction scheme is adopted to uncouple the velocity and the
pressure. We refer to \cite{Guermond2004} for the convergence analysis
of the metod and to \cite{Guermond2005} for a review on projection
methods.

\subsubsection{The Space Discretization}
Let $\bF_D$ be a continuous, piecewise quadratic approximation of
$\bef_D$ on $\Gamma_D$.  The finite element discretization of the
velocity and the pressure is done by using the following linear and
affine spaces:
\begin{align}
  \pmb{\mathbb V}_0(\mathcal{T}) 
&:= \{ \bV \in C^0({\overline{\Lambda}};\mathbb{R}^d) \ | \ \bV{|_{\Gamma_D}} = 0, \ 
\bV\SCAL\bnu |_{\Gamma_S} = 0,\ \bV|_{K} \in \pmb{\mathbb Q}^{2}(K),
  \ \forall K \in \mathcal{T} \},  \\
  \pmb{\mathbb V}_D(\mathcal{T}) 
&:= \{ \bV \in C^0({\overline{\Lambda}};\mathbb{R}^d) \ | \ \bV{|_{\Gamma_D}} = \bF_D, \ 
\bV\SCAL\bnu |_{\Gamma_S} = 0, \ \bV|_{K} \in \pmb{\mathbb Q}^{2}(K),
  \ \forall K \in \mathcal{T} \},  \\
  \mathbb{M}(\mathcal{T})
&:= \{ Q \in C^0({\overline{\Lambda}})
  \rightarrow \mathbb{R} \ | \ Q|_{K} \in \mathbb{Q}^1(K) , \ \forall
  K \in \mathcal{T} \},
\end{align}
Upon setting $\bU(0)=\bU_0$, where $\bU_0$ is a continuous, piecewise
quadratic approximation of the initial velocity $\bu_0$ in
$\pmb{\mathbb V}_D(\mathcal{T})$, the semi-discrete formulation of
\eqref{eqn:navier_stokes} consists of looking for $\bU\in
\mathcal{C}^1([0,T];\pmb{\mathbb V}_D(\mathcal{T}))$ and $P \in
\mathcal{C}^0([0,T];\mathbb{M}(\mathcal{T})$ such that the following
holds for every $t \in (0,T]$:
\begin{equation}\label{eqn:NS_space}
\begin{split}
  \int_{\Lambda}\rho \left(\dfrac{\partial}{\partial t} \bU + (\bU
    \SCAL \nabla)\bU\right) \SCAL \bV + 2\int_{\Lambda} \mu \left(
    \nabla^S \bU{:}\nabla^S\bV \right) -
  \int_{\Lambda}  P \text{div}\bV  + \int_{\Lambda} Q \, \text{div}\bU\\
  = \int_{\Lambda}\rho {\bg}\SCAL \bV + \int_{\Gamma_N} \bef_N \SCAL
  \bV + \int_{\Sigma(t)} \sigma \kappa \, \bn \SCAL \bV, \qquad
  \forall (\bV,Q) \in \mathbb{V}_0(\mathcal{T}) \CROSS
  \mathbb{M}(\mathcal{T}),
\end{split}
\end{equation}

\subsubsection{Time discretization} 
\label{subsec:time_navier_stokes}
The time discretization of \eqref{eqn:NS_space} is done by using the
second-order backward differentiation formula (BDF2) and an
incremental pressure correction scheme in rotational form introduced
and studied in \cite{Guermond2009} to decouple the velocity and the
pressure. In addition to the approximation of the initial condition
${\bU}^0$ mentioned above, the algorithm requires an approximation
$P_0 \in \mathbb M(\mathcal T)$ of the initial pressure $p(0)$.  We
denote by $\bU^n$, $P^n$ and $\Psi^n$ the approximations of
$\bU(\cdot,t^n)$ and $P(\cdot,t^n)$ and the pressure correction,
respectively.

The initialization step consists of setting: $\bU^0=\bU_0$,
$P^{-1}=P^0=P_0$ and $\Psi^0 =0$.  Then, given a new time step $\delta
t^{n+1}$, given $\bU^{n} \in \pmb{\mathbb V}_D(\mathcal T)$, $\Psi^{n}
\in \mathbb M(\mathcal T) $ and $P^{n}\in \mathbb M(\mathcal T)$,
and assuming for the time being that
$\rho(t^{n+1})$, $\mu(t^{n+1})$ and $\Sigma(t^{n+1})$ are also known (see
\S\ref{subsec:comp_alg}), the new fields $\bU^{n+1} \in \pmb{\mathbb
  V}_D(\mathcal T)$, $\Psi^{n+1} \in \mathbb M(\mathcal T)$
and $P^{n+1} \in \mathbb M(\mathcal T)$ are computed in three steps:\\
\step{1}\emph{Velocity Prediction:} Find $\bU^{n+1} \in \pmb{\mathbb
  V}_D(\mathcal T)$ such that
\begin{multline}
\int_{\Lambda} \rho(t^{n+1})\frac {\delta_{\mbox{\tiny BDF2}} }{\delta t^{n+1}} {\bU}^{n+1} \SCAL {\bV}
+ 2\int_{\Lambda}{\mu(t^{n+1})} \left( \nabla^S({\bU}^{n+1}) {:} \nabla^S{\bV} \right) +S_{\mathcal T}(\bU^{n+1},\bV)\\
=-\int_{\Lambda} \rho(t^{n+1}) (({\bU}^n)^{\star} \SCAL \nabla {\bU}^{n}) \SCAL {\bV} 
+\int_{\Lambda}(P^n + \dfrac{4}{3}\Psi^n - \dfrac{1}{3}\Psi^{n-1}) \mbox{div}({\bV})  \\
+ \int_{\Lambda} \rho^{n+1} {\bg^{n+1}}
+ \int_{\Gamma_N} \bef_N^{n+1} \SCAL \bV 
+\int_{\Sigma(t^{n+1})} \sigma \kappa^{n+1} \, \bn^{n+1} \SCAL \bV,  \quad
\forall \bV \in \pmb{\mathbb V}_0(\mathcal{T}),
\label{eqn:space:first}
\end{multline}
where 
\begin{itemize}
\item the BDF2 approximation of the time derivative with variable time
  stepping is given by
\begin{equation*}
  \dfrac{\partial}{\partial t}\bU(\bx,t^{n+1}) \approx  \frac {\delta_{\mbox{\tiny BDF2}} }{\delta t^{n+1}}  
  \bU^{n+1}(\bx):= \dfrac{1}{\delta t^{n+1}} 
  \left( \dfrac{1+2 \eta_{n+1}}{1 + \eta_{n+1}} \bU^{n+1}(\bx) - (1+ \eta_{n+1})\bU^{n}(\bx) + \dfrac{\eta^2_{n+1}}{1 + \eta_{n+1}} \bU^{n-1}(\bx)\right), \qquad \eta_{l+1} := \frac{\delta t^{l+1}}{\delta t^l};
\end{equation*}
\item the extrapolated velocity $(\bU^n)^\star$ is defined by
$(\bU^n)^{\star} := \bU^{n} + \eta^{n}( \bU^n - \bU^{n-1})$;
\item as discussed in \cite{enumath}, the bilinear form $S_{\mathcal
    T}:\mathbb V(\mathcal T)\CROSS \mathbb V(\mathcal T)\rightarrow
  \mathbb R$ is added to control the divergence of the velocity and to
  cope with variable time stepping and open boundary conditions
\begin{equation}\label{eqn:stab}
S_{\mathcal T}(\bW,\bV) := 
C_\text{stab} \sum_{K\in\mathcal{T}}\int_{K}  (\mu^{n+1} + \rho^{n+1}  \|h_K \bW \|_{L^{\infty}({K})}) 
(\nabla \SCAL \bW) (\nabla \SCAL \bV),
\end{equation}
where $C_\text{stab}$ is an absolute constant.
\end{itemize}

\step{2}\emph{Pressure Correction Step:} The pressure increment
$\Psi^{n+1} \in \mathbb{M}(\mathcal{T})$ is determined by solving
\begin{equation}
\int_{\Lambda} \nabla \Psi^{n+1} \nabla Q
= -\dfrac{3 \min_{x\in \Lambda} \rho(t^{n+1})}{2 \delta t^{n+1}}  \int_{\Lambda} \text{div}({\bU}^{n+1}) Q ,  \ \ \ \forall Q \in \mathbb{M}(\mathcal{T}).
\label{eqn:space:second}
\end{equation}

\step{3}\emph{Pressure Update:} 
The pressure $P^{n+1} \in \mathbb{M}(\mathcal{T})$ is obtained by solving 
\begin{equation}
\int_{\Lambda}P^{n+1} Q 
  = \int_{\Lambda} \left( P^n  + \Psi^{n+1} \right)Q 
  - \min_{x\in \Lambda}\mu(t^{n+1}) \int_\Lambda \mbox{div}({\bU}^{n+1}) \, Q ,  \ \forall Q \in \mathbb{M}(\mathcal{T})  .
\label{eqn:space:third}
\end{equation}

For stability purposes, we restrict the space and time discretization parameters to satisfy a CFL condition
\begin{equation}\label{e:CFL2}
\delta t^{n+1} \leq C_{\textrm{CFL}} \frac 1 2 \frac{\min_{K\in \mathcal T_h} h_K}{ \| \bU^{n}\|_{L^\infty(\Lambda)}},
\end{equation}
where $C_{\textrm{CFL}}$ is the same constant appearing in
\eqref{eqn:CFL} and $\frac 1 2 \min_{K\in \mathcal T_h} h_K$
is the minimum distance between two Lagrange nodes using
$\polQ_2$ elements.

\subsection{The Surface Tension}
\label{subsec:surf_tens}
This section describes the approximation of the curvature term
$\int_{\Sigma(t^{n+1})} \sigma \kappa^{n+1} \, (\bn^{n+1} \SCAL \bV)$
appearing in the first step of the projection
method \eqref{eqn:space:first}.  We follow the method proposed by \cite{Bansch} (see also
\cite{hysing2009quantitative}) using the work of \cite{Dziuk}.
This approach is based on the following representation of total curvature:
\[
\kappa \bn = \nabla_\Sigma \SCAL \nabla_\Sigma  \mathrm{\bf Id}_\Sigma,
\]
where $\mathrm{\bf Id}_\Sigma$ is the identity mapping on $\Sigma$
and, given any extension $\tilde \bv$ of $\bv$ in a neighborhood of
$\Sigma$, the tangential gradient of $\bv:\Sigma \rightarrow \mathbb
R^d$ is defined by $\nabla_\Sigma \bv := \nabla {\tilde{\bv}}|_\Sigma
( \text{\bI} - {\bn} \otimes {\bn})$, see \eg \cite{Gilbarg}.
Multipliying the above identity by a test function $\bV$ and
integrating by parts over $\Sigma$ yields
\begin{equation}
\int_{\Sigma} \sigma \kappa {\bn} \SCAL {\bV}  
= -\int_{\Sigma} \sigma {\nabla}_{\Sigma}\mathrm{\bf Id}_\Sigma : {\nabla}_{\Sigma} {\bV} 
   +\int_{\partial \Sigma} \sigma \partial_\btau \mathrm{\bf Id}_\Sigma \SCAL \bV,
\label{eqn:stokes:st_1}
\end{equation}
where $\btau$ is the co-normal to $\Sigma$ and $\partial_{\btau}$ is
the derivative in the co-normal direction.  In the present context the
integral on $\partial \Sigma$ vanishes since either $\Sigma$ is a
closed manifold or $\partial_\btau \mathrm{\bf Id}_\Sigma=0$.
This identity together with the first-order prediction $\mathrm{\bf
  Id}_{\Sigma(t^{n+1})} \approx \mathrm{\bf Id}_{\Sigma(t^n)} + \delta
t^{n+1} \bU^{n+1}$ of the interface evolution
\eqref{eqn:free_surface} gives a semi-implicit representation of the
total curvature
\begin{align*}
  \int_{\Sigma(t^{n+1})} \sigma \kappa^{n+1} \, (\bn^{n+1}\SCAL\bV)
  &= - \int_{\Sigma(t^{n+1})}  {\nabla}_{\Sigma(t^{n+1})}\mathrm{\bf Id}_{\Sigma(t^{n+1})}{:}{\nabla}_{\Sigma(t^{n+1})} {\bV}  \\
  & \approx -\int_{\Sigma(t^{n+1})} {\nabla}_{\Sigma(t^{n+1})}\mathrm{\bf
    Id}_{\Sigma(t^{n+1})} {:} {\nabla}_{\Sigma(t^{n+1})} {\bV} - \delta
  t^{n+1} \int_{\Sigma(t^{n+1})} \sigma
  {\nabla}_{\Sigma(t^{n+1})}\bU^{n+1}{:}{\nabla}_{\Sigma(t^{n+1})} {\bV}.
\end{align*}
One key benefit of this representation comes from the additional
stabilizing term $\int_{\Sigma(t^{n+1})} \sigma
{\nabla}_{\Sigma(t^{n+1})}\bU^{n+1} {:} {\nabla}_{\Sigma(t^{n+1})}
{\bV}$, which we keep implicit in \eqref{eqn:space:first}. The
technicalities regarding the approximation of this integral using the
level set representation are detailed in section \ref{Sec:data_NS},
see \eqref{discrete_surface_tension}.

\subsection{The Coupled System}
\label{subsec:comp_alg}
\begin{figure}[ht]
\centering
\begin{tikzpicture}[node distance = 3cm, auto]
    \node [cloud] (Init) {\footnotesize Initialization};
    \node [block, right of=Init] (sol_level) {\footnotesize Solve Level-set};
    \node [block, right of=sol_level] (sol_vel) {\footnotesize Solve Velocity \\ Update Pressure};
    \node [block, right of=sol_vel] (end) {\footnotesize Refine Mesh \\ Choose Time Step};

    \path [line] (Init) -- (sol_level);
    \path [line] (sol_level) -- (sol_vel);
    \path [line] (sol_vel) -- (end);
    \draw [-] (10,-0.52) -- (10,-0.75);
    \draw (10,-0.75) -- (3.35,-0.75);
    \path [line] (3.35,-0.75) -- (3.35,-0.55);
    \draw [blue] (1.5,-1) rectangle (11.5,0.75);
\end{tikzpicture}
\caption{Global Algorithm Flowchart.}
\label{fig:arg_flow}
\end{figure}
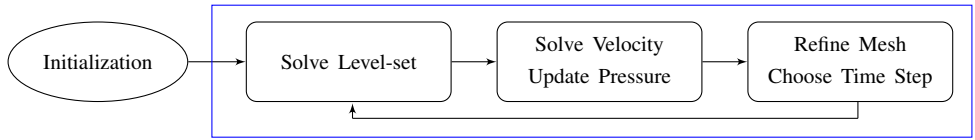

The two solvers described in Sections \ref{subsec:num_level_set} and
\ref{ss:NS} above are sequentially combined.  The flowchart of the
resulting free boundary flow solver is shown in Figure
\ref{fig:arg_flow}.  The remaining subsections of
\S\ref{subsec:comp_alg} detail the coupling between the two solvers
and other implementation technicalities.

\subsubsection{Data for the Level-Set Solver \eqref{eqn:level_ssprk3}-\eqref{eqn:level_ssprk3_2}}
Following \eqref{eqn:level_ssprk3}, given $\Phi^n$, an approximation
of $\phi(t^n)$, the approximation $\Phi^{n+1}$ of $\phi(t^{n+1})$
requires the velocities $\bu(t^n), \bu(t^{n+\frac{1}{2}})$ and
$\bu(t^{n+1})$.  To avoid an implicit coupling between the level-set
solver and the Navier-Stokes solver, these quantities are replaced by second-order
extrapolations using $\bU^n$ and $\bU^{n-1}$:
$$
\bu(t^n) \approx \bU^n, \qquad \bu(t^{n+\frac12})\approx \bU^{n}+\frac 1 2\frac{\delta t^{n+1}}{\delta t^n} (\bU^n-\bU^{n-1}), \qquad \text{and} \qquad  \bU^{n+1} \approx \bU^{n}+\frac{\delta t^{n+1}}{\delta t^n} (\bU^n-\bU^{n-1}).
$$

The $\mbox{sign}(.)$ function that is used in the right-hand side of
\eqref{eqn:level:final} to make sure that $\|\nabla\phi\|_{\ell^2}$ is
close to 1 in a small neighborhood of $\Sigma$ (\ie on the fly
reinitialization, see discussion following \eqref{eqn:level:sign}) is
redefined and replaced by:
\begin{equation}\label{eqn:level:sign:apprx}
\mbox{sign}_h(s) = \left\lbrace 
\begin{array}{ll}
+1 & \qquad \text{if }s> \beta \tanh(C_S), \\
-1 & \qquad \text{if }s< -\beta \tanh(C_S), \\
0 & \qquad \text{otherwise,}
\end{array}
\right. 
\end{equation}
where $C_S$ is an absolute constant.  The thresholding in the above
definition of the approximate $\mbox{sign}$ function gives
$\mbox{sign}_h(\phi(\bx,t))=\pm 1$ whenever
$\frac{|\phi(\bx,t)|}{\beta} \ge \tanh(C_S)$,
which is compatible with the behavior $\frac{\phi(\bx,t)}{\beta}
\approx \tanh\left(\frac{\text{dist}(\Sigma(t),\bx)}{\beta}\right)$
that is expected for the level-set function, see
\eqref{eqn:cut_off_beta}.

The parameter $\lambda$ (``reinitialization relative speed'' in the
language of \cite{ville2011convected}) is defined for $t \in
(t^n,t^{n+1}]$ by 
\begin{equation}\label{eqn:level:lambda} 
\lambda = C_\lambda \| \bU^n\|_{L^\infty(\Lambda)}, \qquad K\in \mathcal T,
\end{equation}
where $C_\lambda$ is an absolute constant.  This definition is
motivated by the CFL condition \eqref{eqn:CFL}.  We refer to
\cite{ville2011convected} for further details.

\subsubsection{Data for the Navier-Stokes Solver}
\label{Sec:data_NS}
The definitions of the fields $\bU^{n+1}$ and $P^{n+1}$ in
\eqref{eqn:space:first}--\eqref{eqn:space:third} invoke the values of
the density field $\rho^{n+1}$ and viscosity field $\mu^{n+1}$.  Once
$\phi^{n+1}$ is computed, these quantities are evaluated by using the following definitions
\begin{align}
 \rho^{n+1} &= \rho^+ \dfrac{1+H_h(\Phi^{n+1})}{2} + \rho^{-} \dfrac{1-H_h(\Phi^{n+1})}{2} \\
 \mu^{n+1} &= \mu^+ \dfrac{1+H_h(\Phi^{n+1})}{2} + \mu^{-} \dfrac{1-H_h(\Phi^{n+1})}{2}  ,
\label{eqn:rho_mu}
\end{align}
where, $\rho^{\pm}, \mu^{\pm}$ are the density/viscosity in
$\Omega^{\pm}$, and $H_h(.)$ an approximation of the Heaviside
function defined as follows:

\begin{equation}
H_h(s)=\left\{ 
{\begin{array}{cl}
1,  & \qquad \text{if } s > \beta \tanh(C_H), \\
-1, & \qquad \text{if } s < -\beta \tanh(C_H), \\
\dfrac{s}{\beta \tanh(C_H)}, & \qquad \text{otherwise},
\end{array}} 
\right.
\label{eqn:H_level}
\end{equation}
where $C_H$ is an absolute constant, as suggested in
\cite{ville2011convected}.  Similarly to what we have done to
approximate the $\mbox{sign}$ function, the above regularization is
compatible with the behavior $\frac{\phi(\bx,t)}{\beta} \approx
\tanh\left(\frac{\text{dist}(\Sigma(t),\bx)}{\beta}\right)$ that is
expected for the level-set function.

The approximation of the surface tension term in
\eqref{eqn:space:first} is done by following
\cite{hysing2009quantitative}.  Let $\epsilon>0$, and consider the
piecewise linear regularized Dirac measure supported on
$\Sigma(t^{n+1})$, say $\delta_{\epsilon}$, defined by
\begin{equation}
\delta_{\epsilon}(\bx) := 
\left\lbrace 
\begin{array}{ll}
\dfrac{1}{\epsilon}\left( 1 - \frac{\text{dist}(\Sigma,\bx)}{\epsilon} \right) & \qquad \text{if}\quad  |\text{dist}(\Sigma,\bx)| < \epsilon,   \\
 0                                      & \qquad \text{otherwise}.
\end{array} \right.
\label{eqn:st:delta}
\end{equation}
To account for the fact that we do not have access to the distance
to the interface $\Sigma$ but rather to an approximation of
$\beta \tanh(\frac{\text{dist}(\Sigma,\bx)}{\beta})$, see Section
\ref{subsec:levelset:reinit}, we rescale $\delta_{\epsilon}(\bx)$ as suggested in
\cite{Engquist04discretizationof,tornberg2000interface} and consider instead:
\begin{equation}
  \delta_{\epsilon}(\bx,\Phi) = 
  \left\{ \begin{array}{ll}
      \dfrac{1}{\tilde \epsilon}\left( 1 - \frac{\Phi(\bx)}{\tilde \epsilon} \right) 
      \|\nabla \Phi(\bx)\|_{\ell^2},  
      &|\Phi(\bx)| < \tilde{\epsilon} := \epsilon \frac{\|\nabla \Phi
        (\bx)\|_{\ell^1}}{\| \nabla \Phi (\bx)\|_{\ell^2}}\\
      0,                                        & \text{otherwise,}
\end{array} \right.
\end{equation}
where $\|\cdot\|_{\ell^p}$ is the $\ell^p$-norm in $\mathbb R^d$.  In
practice, we chose $\epsilon=\beta \tanh(C_H)$ to be consistent with
the approximation of the Heaviside function, see \eqref{eqn:H_level}.
Using this approximate Dirac measure, the approximation of the surface
tension discussed in Section \ref{subsec:surf_tens} becomes
\begin{align}
  \int_{\Sigma(t^{n+1})} \sigma \kappa^{n+1} \, (\bn^{n+1}\SCAL\bV)
  & \approx -\int_{\Sigma(t^{n+1})} {\nabla}_{\Sigma(t^{n+1})}\mathrm{\bf
    Id}_{\Sigma(t^{n+1})} {:} {\nabla}_{\Sigma(t^{n+1})} {\bV} - \delta
  t^{n+1} \int_{\Sigma(t^{n+1})} \sigma
  {\nabla}_{\Sigma(t^{n+1})}\bU^{n+1}{:}{\nabla}_{\Sigma(t^{n+1})} {\bV} \nonumber \\
  & \approx -\int_{\Lambda}\left( {\nabla}_{\Phi^{n+1}}\mathrm{\bf
    Id}_{\Lambda} {:} {\nabla}_{\Phi^{n+1}} {\bV}\right) \delta_\epsilon(.,\Phi^{n+1})- \delta
  t^{n+1} \int_{\Lambda} \left(\sigma
  {\nabla}_{\Phi^{n+1}}\bU^{n+1}{:}{\nabla}_{\Phi^{n+1}} {\bV}\right) \delta_\epsilon(.,\Phi^{n+1}), \label{discrete_surface_tension}
\end{align}
where $\mathrm{\bf
    Id}_{\Lambda}$ is the identity mapping on $\Lambda$ and 
\[
{\nabla}_{\Phi^{n+1}} \mathrm{\bf  Id}_{\Lambda} 
= I-\frac{\nabla \Phi^{n+1} \otimes \nabla \Phi^{n+1}}{\| \nabla \Phi^{n+1}\|_{\ell^2}^2},\qquad
{\nabla}_{\Phi^{n+1}}\bV 
:= \nabla \bV \left(I-\frac{\nabla \Phi^{n+1} \otimes \nabla \Phi^{n+1}}{\| \nabla \Phi^{n+1}\|_{\ell^2}^2}\right),
\quad \forall \bV  \in \pmb{\mathbb V}_D(\mathcal T).
\]
Note that the above definitions correspond to approximating the normal
vector $\bn$ on
$\Sigma(t^{n+1})$ by $-\frac{\nabla \Phi^{n+1}}{\| \nabla \Phi^{n+1}\|_{\ell^2}}$.


\subsection{Adaptive Mesh Refinement}
\label{subsec:adaptivity}
The experiments reported in \cite{BounceJet,B_2009,PhysRevE.87.061001}
suggest that a critical feature of bouncing jets is the occurrence of
a thin layer between the jet and the rest of the fluid. These
observations have led us to adopt a mesh refinement technique to
describe accurately this thin layer.  The refinement strategy is
designed to increase the mesh resolution around the zero-level set of
$\Phi$.  More precisely, a cell $K$ is refined if its generation count
(the number of times a cell from the initial subdivision has been
refined to produce the current cell) is smaller than a given number
$R_{\max}$ and if
\begin{equation}
|\Phi(\bx_K,t)| \leq  \beta \tanh(C_R),
\label{eqn:level:x_r}
\end{equation}
where $\bx_K$ is the barycenter of $K$ and $C_R$ is an absolute
constant.  The purpose of the parameter $R_{\max}$ is to control the
total number of cells. The parameter $C_R$ controls the distance to
the interface below which refinement occurs.  Note that
\eqref{eqn:level:x_r} is compatible with \eqref{eqn:level:sign:apprx}
and \eqref{eqn:H_level}.  The subdivisions are done with at most one
hanging node per face, see \eg \cite[Section~6.3]{MR2670003}.  A cell
is coarsen if it satisfies the following three conditions: its
generation count is positive;
\begin{equation}
|\Phi(\bx_K,t)| \geq  \beta \tanh(C_C),
\label{eqn:level:x_c}
\end{equation}
where $C_{C}\geq C_R$ is an absolute constant; and if once coarsened,
the resulting subdivision does not have more than one hanging node per
face.  We refer to the documentation of the deal.II library for
further details, \cite{BangerthHartmannKanschat2007}.

\section{Numerical Validations}\label{s:validation}
The algorithm presented in the previous sections has been implemented
using the \emph{deal.II} finite element library described in
\cite{BangerthHartmannKanschat2007,dealII81}. Parallelism is handled
by using the MPI (Message Passing Interface) library, see
\cite{gabriel04:_open_mpi}. The subdivision and mesh distribution is
done by using the p4est library from \cite{p4est}.

The rest of this section illustrates and evaluates the performance of
the above algorithm.  We start in Section \ref{subsec:constants} by
specifying all the numerical constants required by the algorithm.  The
validation of the transport code for solving the level-set equation in
done Section~\ref{sec:validation_transport}. The validation of the
two-phase fluid system is done in Section~\ref{sec:validation_two}.

\subsection{Numerical Parameters}\
\label{subsec:constants}
Our algorithm involves several numerical parameters.  In this section,
we briefly recall their meaning, where they appear, and we specify the
value of each of them. Unless specified otherwise these values are
fixed for this entire section.
\begin{table}[ht]
\footnotesize
\centering
\begin{tabular}{|r || c |  c |  c  | c |  c | c | c | c | c | c |c|}
\hline
&$C_{\text{CFL}}$ & $C_\lambda$ &  $C_H$ & $C_S$ &$C_R $ & $C_C$ & $R_{\max}$ & $C_\text{Lin}$ & $C_\text{Ent}$& $p$  &$C_{stab}$ \\  \hline\hline
Purpose & CFL & \multicolumn{3}{c|}{Reinitialization} &  \multicolumn{3}{c|}{Adaptivity} &  \multicolumn{4}{c|}{Stabilization} \\ \hline
Appears in & \eqref{eqn:CFL}\&\eqref{e:CFL2}   & \eqref{eqn:level:lambda} & \eqref{eqn:H_level} &  \eqref{eqn:level:sign:apprx}
 &  \eqref{eqn:level:x_r} & \eqref{eqn:level:x_c} &\S\ref{subsec:adaptivity} &
 \eqref{eqn:lin:visc} & \eqref{eqn:level:fv}&  \eqref{eqn:entropy:choice} & \eqref{eqn:stab} \\ \hline
Value  & 0.25  & 0.01       &   1.25   & 0.5 &   2. & 2. & 2. & 0.1 & 0.1& 20  & 0.1 \\  \hline
\end{tabular}
\caption{Description and values of the numerical parameters. }
\label{tab:num_param}
\end{table}
A complete list of all the parameter is shown in
Table~\ref{tab:num_param}.  The parameter $\beta$ determining the
width of the hyperbolic tangent filter \eqref{eqn:cut_off_beta} is
defined to be $\beta= \min_{K\in\mathcal{T}}h_K$.  This definition is
used in the expression of the hyperbolic tangent filter
\eqref{eqn:cut_off_beta}, in the threshold of the approximate Heaviside
function  \eqref{eqn:H_level} and the
approximate \mbox{sign} function \eqref{eqn:level:sign:apprx}, and in the
refinement strategy \eqref{eqn:level:x_r}-\eqref{eqn:level:x_c}.

\subsection{Transport of the Level-Set} \label{sec:validation_transport}
\subsubsection{Convergence Tests}
The consistency of the algorithm for the approximation of the
level-set equation \eqref{eqn:level_ssprk3}-\eqref{eqn:level_ssprk3_2}
is a priori third-order in time and space.  To evaluate whether this
is indeed the case, we solve the linear transport equation in the unit
square $\Lambda=(0,1)^2$ using the velocity field 
\begin{equation*}
\bu(x,y,t) := \left( 
\begin{array}{r}
-\sin^2(\pi x)\sin(2\pi y)\cos(\pi t/0.1) \\
\sin^2(\pi y)\sin(2\pi x)\cos(\pi t/0.1)
\end{array}\right).
\end{equation*}
This flow is time-periodic of period 1, wich implies that
$\phi(\bx,0.2)=\phi_0(\bx)$, for all $\bx\in \Lambda$.  The initial
level-set function is chosen to be the signed distance to the line of
equation $y=0.5$:
\begin{equation}
\phi_0(x,y) =  y - 0.5.
\label{eqn:level:init_test}
\end{equation}
The errors are evaluated at $t=0.2$.
Three different scenarios are considered: (i) no stabilization and no
reinitialization; (ii) entropy viscosity stabilization and no
reinitialization; (iii) entropy viscosity and reinitialization, \ie
the complete algorithm. Except for $\beta$,
the values of the numerical parameters are given in Table
\ref{tab:num_param}.
%
We consider four computations done on four uniform meshes with
constant time steps.  The mesh-size and the time step are divided by
$2$ each time. The meshes are composed of 1089, 4225, 16641and 66049
$\polQ_2$ degree of freedoms.  The space discretization is chosen fine
enough not to influence the time error. We report in Table
\ref{table:levelset_vortex}, the errors
$\|\Phi^N-\phi_0\|_{L^2(\Lambda)}$ for the three scenarios and the
observed rates of convergence.  Note that the reinitialization is
turned on in scenario (iii), which implies that the exact solution at
$t=0.2$ is not $\phi_0$ anymore; we must instead compare $\Phi^N$ with
$\beta \tanh\left(\frac{\phi_0(x,y)}{\beta}\right)$. In this
particular case we keep $\beta$ constant to ascertain the third-order
consistency of the algorithm; we set $\beta = 0.0203125$.
\begin{table}[!h]
\centering
\begin{tabular}{|c||c|c||c|c||c|c|} \hline
& \multicolumn{2}{c||}{Scenario (i)}
& \multicolumn{2}{c||}{Scenario (ii)}
& \multicolumn{2}{c|}{Scenario (iii)} \\
\hline
$\delta t$   & $L^2$ error &  Observed rate 
& $L^2$ error &  Observed rate
& $L^2$ error &  Observed rate \\ \hline
1e-2                   & 1.2e-5     &   -   
                       & 1.5e-5    &  - 
                       & 1.3e-3    &  -\\ \hline
5e-3                   & 1.5e-6     &  3.
                       & 1.8e-6     &  3.
                       & 1.5e-4     &  3.\\ \hline
2.5e-3                 & 1.9e-7     & 3.
                       & 2.0e-7    &  3.2
                       & 2.0e-5    & 3.\\ \hline
1.25e-3                & 2.4e-8     &   3.
                       & 2.5e-8     &  3.1
                       & 2.6e-6    &   3. \\ \hline
\end{tabular}
\caption{Single vortex test case: $L^2$-norm of error and observed convergence rates 
for each scenario after one period for different time discretization resolutions. 
(i) RK3 Galerkin algorithm; (ii) RK3 algorithm with entropy viscosity; 
(iii) RK3 with entropy viscosity and reinitialization.
  Third-order convergence rate is observed in all scenarios. 
}
\label{table:levelset_vortex}
\end{table}

\subsubsection{Rotation of a circular level-set: Effect of Different
  Viscosities and Reinitialization}
\begin{figure}[ht]
\centering
\begin{subfigure}[]{0.18\textwidth}
\includegraphics[width=\textwidth]{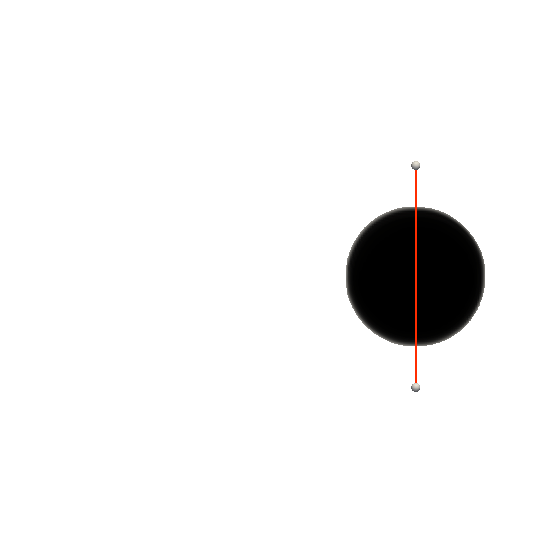}
\caption{}
\label{fig:circle:stab_a}
\end{subfigure}
\begin{subfigure}[]{0.15\textwidth}
\includegraphics[width=\textwidth]{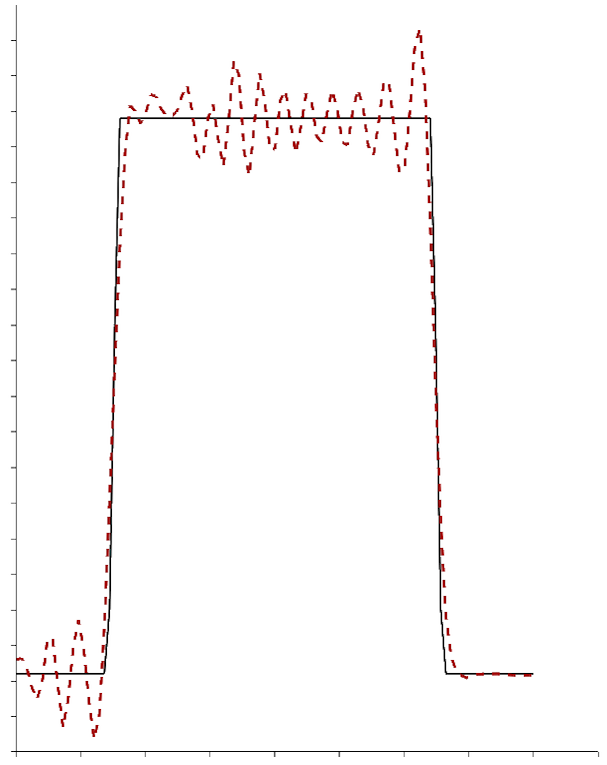}
\caption{}
\label{fig:circle:stab_b}
\end{subfigure}
\begin{subfigure}[]{0.15\textwidth}
\includegraphics[width=\textwidth]{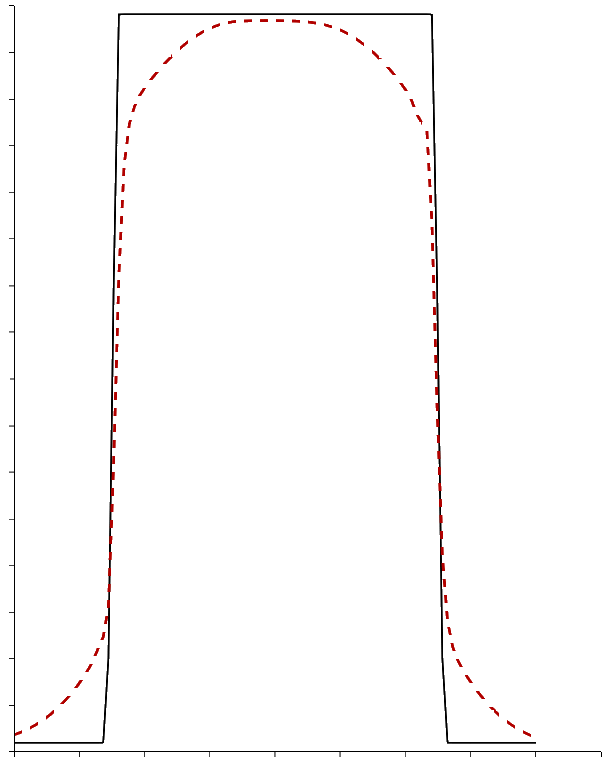}
\caption{}
\label{fig:circle:stab_c}
\end{subfigure}
\begin{subfigure}[]{0.15\textwidth}
\includegraphics[width=\textwidth]{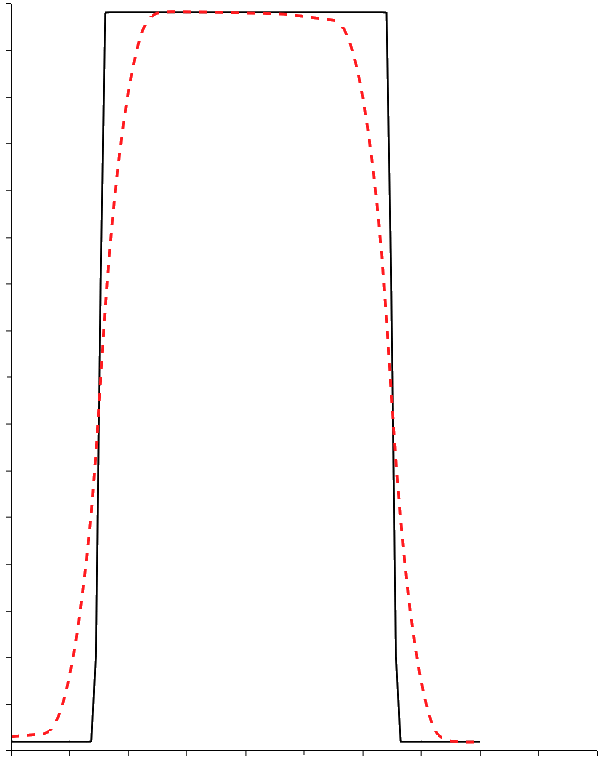}
\caption{}
\label{fig:circle:stab_d}
\end{subfigure}
\begin{subfigure}[]{0.15\textwidth}
\includegraphics[width=\textwidth]{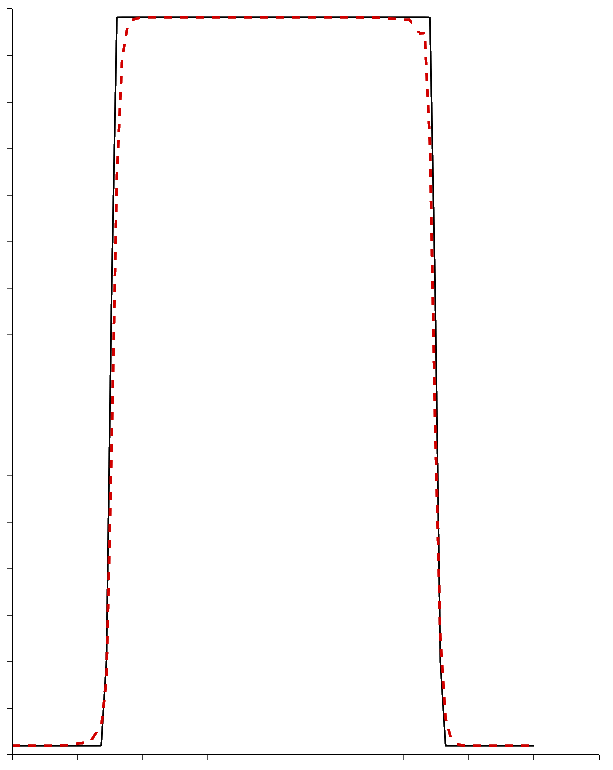}
\caption{}
\label{fig:circle:stab_e}
\end{subfigure}
\caption{(a) Initial level-set: the dark region corresponds to
  $\phi_0\geq 0$; (b)-(e) Comparison of the graphs of the exact (solid
  line) and approximate (dashed line) level-set functions along the
  line of equation $x=0.5$ after one revolution: (b) without
  stabilization and reinitialization
  ($C_{\text{Lin}}=C_{\text{Ent}}=C_\lambda=0$); (c) with first-order
  linear stabilization, $C_{\text{Lin}}=0.1$; (d) With entropy
  residual stabilization, $C_{\text{Lin}}=0.1$ and
  $C_{\text{Ent}}=0.1$; (e) with entropy residual stabilization and
  reinitialization, $C_\lambda=0.01$. The transport without
  stabilization of the nearly discontinuous level-set function yields
  spurious oscillations.  These oscillations are removed by the linear
  viscosity at the expense of large numerical diffusion.  The
  numerical diffusion is minimized by turning on the entropy residual
  term.  Finally, adding the reinitialization allows to nearly recover
  the exact profile.  }
\label{fig:circle:stab}
\end{figure}

In this test case the initial data for the level set is the hyperbolic
tangent filter applied to the distance to the circle centered at
$(0.5,0)$ and of radius $0.25$
\[
\phi_0(x,y) := \beta \tanh\left(\frac{0.25-((x-0.5)^2 + y^2)^\frac12}{\beta} \right),
\qquad (x,y) \in \Lambda:=(-1,1)^2,
\]
see Figure \ref{fig:circle:stab_a}. The level-set is transported by
using the solid rotation velocity field:
\begin{equation}
\bu(x,y) := \left( 
\begin{array}{r}
-y \\
x
\end{array}\right).
\label{eqn:level:circle_vel}
\end{equation}

The computation is done with the
adaptive algorithm described in Section \ref{subsec:adaptivity} using
$R_{\max}=2$.  The initial mesh is uniform with $h_K = 0.015625$ for
all $K\in \mathcal T$. The resulting time-dependent mesh is such that
$\min_{T\in \mathcal T} h_K \approx 0.00390625$, \ie we set $\beta =
0.00390625$.  The values of the other numerical constants are provided
in Table \ref{tab:num_param}. The time step $\delta t$ is chosen to be
uniform and equal to $\pi\CROSS 10^{-3}$.  We compare the initial
level-set with its approximation after one revolution.
Figures~\ref{fig:circle:stab_b}-\ref{fig:circle:stab_e} illustrate the
benefits of using the entropy stabilization and the reinitialization
technique by comparing the graphs of the exact and approximate
level-set functions along the line of equation $x=0.5$ after one
revolution.

\subsubsection{3D Slotted Disk: Long Time Behavior}
\begin{figure}[!b]
\centering
\begin{subfigure}[]{0.146\textwidth}
\includegraphics[width=\textwidth]{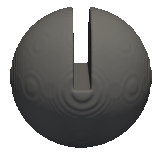}
\caption{}
\label{fig:zalesak_a}
\end{subfigure}
\begin{subfigure}[]{0.15\textwidth}
\includegraphics[width=\textwidth]{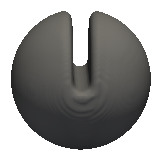}
\caption{}
\end{subfigure}
\begin{subfigure}[]{0.15\textwidth}
\includegraphics[width=\textwidth]{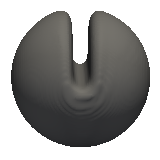}
\caption{}
\end{subfigure}
\begin{subfigure}[]{0.15\textwidth}
\includegraphics[width=\textwidth]{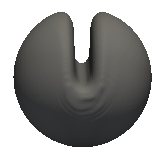}
\caption{}
\end{subfigure}
\begin{subfigure}[]{0.15\textwidth}
\includegraphics[width=\textwidth]{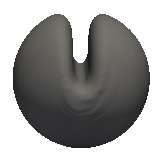}
\caption{}
\end{subfigure}
\begin{subfigure}[]{0.15\textwidth}
\includegraphics[width=\textwidth]{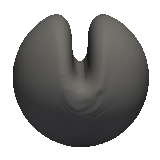}
\caption{}
\end{subfigure}
\caption{Rotating Slotted three dimensional sphere. From left to right the dark regions corresponds to $\Phi\geq 0$ after 0,1,2,3,4 and 5 full rotations.
Oscillations and numerical diffusion are controlled by the entropy viscosity as well as the reinitialization algorithm.}
\label{fig:z_disk_3D}
\end{figure}

\begin{figure}[!ht]
\centering
\includegraphics[scale = 0.17]{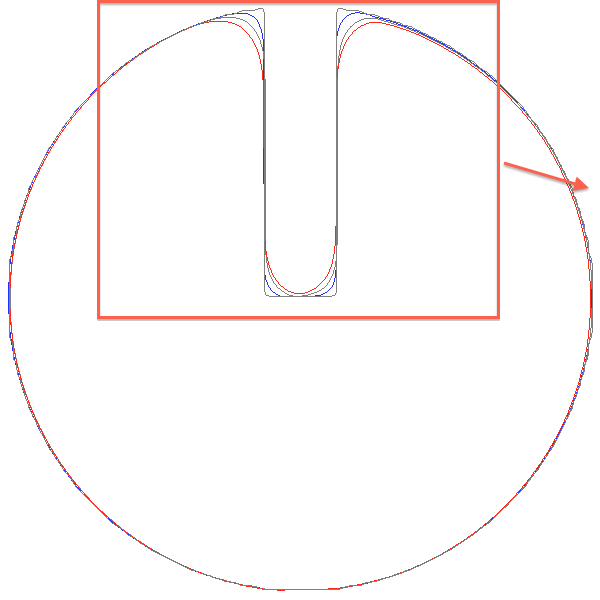} \hspace*{0.1in}
\includegraphics[scale = 0.17]{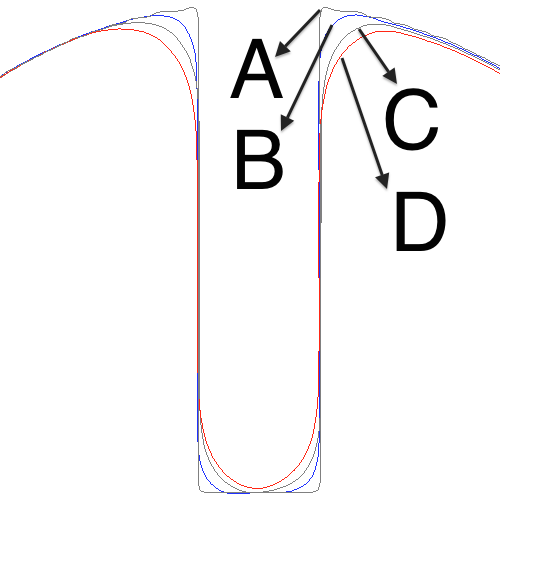}
\caption{Iso-contour $\Phi=0$ in the plane $z=0$: (A) initial data; (B)
  after one rotation; (C) after three rotations; (D) after five
  rotations.  Numerical diffusion is observed but is greatly minimized
  by the entropy viscosity and reinitialization techniques.}
\label{fig:z_disk:cont}
\end{figure}

A typical benchmark for the transport of a level-set function is the
so-called Zalesak disk documented in \citep{Zalesak1979}. We consider
in this section the three-dimensional version thereof. The
computational domain is $\Lambda:=(0,1)^3$.  The initial level-set is
the characteristic function of a slotted sphere centered at
$(0.5,0.75,0.5)$ with a radius of $0.15$.  The width, height and depth
of the slot are $0.0375$, $0.15$, $0.3$ respectively; see Figure
\ref{fig:zalesak_a}).  The initial profile is transported by using the
following velocity field:
\begin{equation}
\bu(x,y,z) := \left( 
\begin{array}{r}
0.5-y \\
x -0.5 \\
0
\end{array}\right) .
\end{equation}
The time step is chosen to be uniform $\delta t= \pi\CROSS 10^{-3}$.
The initial subdivision is composed of cells of diameter $0.015625$.
The adaptive mesh refinement technique described in Section
\ref{subsec:adaptivity} is used with $R_{\max}=2$; the minimum mesh
size is $\min_{K\in \mathcal T} h_K = 0.00390625$.  The numerical
constants are given in Table \ref{tab:num_param}.
The computation is done until the slotted sphere has undergone 5 full
revolutions. The iso-surface $\Phi =0$ is shown in Figure
\ref{fig:z_disk_3D} after each of the 5 periods.  Oscillations and
numerical diffusion are controlled by the entropy viscosity and the
reinitialization algorithm.  A closer look at the slotted region is
provided in Figure \ref{fig:z_disk:cont}.

\subsubsection{Single Vortex: Large Deformations} 
\begin{figure}[!h]
\centering
\begin{subfigure}[t=0]{0.19\textwidth}
\includegraphics[width=\textwidth,bb=50 65 280 250,clip=]{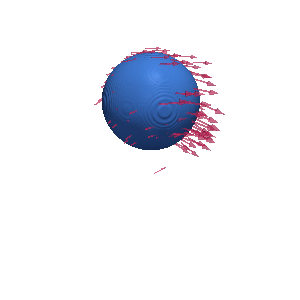}
\end{subfigure}
\begin{subfigure}[t=1]{0.19\textwidth}
\includegraphics[width=\textwidth,bb=50 65 280 250,clip=]{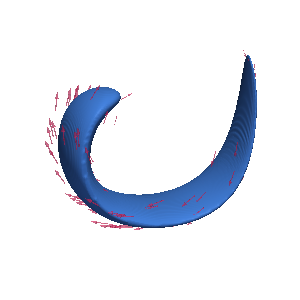}
\end{subfigure}
\begin{subfigure}[t=3]{0.19\textwidth}
\includegraphics[width=\textwidth,bb=50 65 280 250,clip=]{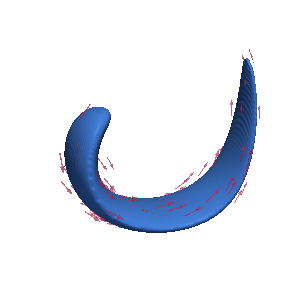}
\end{subfigure}
\begin{subfigure}[t=4]{0.19\textwidth}
\includegraphics[width=\textwidth,bb=50 65 280 250,clip=]{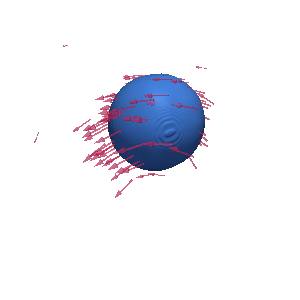}
\end{subfigure}
\caption{3D Vortex using \eqref{eqn:level:vortex_vel}.  From left to
  right: $t=0$, $t=1$, $t=3$, $t=4$.  The initial sphere is recovered after one
  cycle.  The vectors indicate the direction of the velocity field.}
\label{fig:vortex_3d}
\end{figure}

The Single Vortex problem consists of the deformation of a sphere by a
time-periodic incompressible vortex-like flow. 
The computational domain is $\Lambda=(0,1)^3$, and the
time-periodic velocity field is defined by
\begin{equation}
\bu(x,y,t) := \left( 
\begin{array}{c}
-\sin^2(\pi x)\sin(2\pi y)\cos(\pi t/4) \\
\sin^2(\pi y)\sin(2\pi x)\cos(\pi t/4) \\
0
\end{array}\right).
\label{eqn:level:vortex_vel}
\end{equation}
The initial level-set is given by
\begin{equation}
\phi_0(\bx):= \beta \tanh\left(\frac{\mbox{dist}(S,\bx)}{\beta}\right),
\end{equation}
where $S$ is the sphere centered at $(0.5,0.75,0.5)$ of radius of
$0.15$. The field $\phi_0$ is a regularized version of the distance
function using the $\tanh$ cut-off filter \eqref{eqn:cut_off_beta}.
The divergence-free velocity field severely deforms the level-set
until $t=2$ and returns it to its initial shape at $t=4$.
The time step is chosen to be uniform $\delta t=0.001$ and the final
time is $t=4$ (1 cycle).  The initial subdivision is made of uniform
cells of diameter $0.015625$; the minimum mesh size allowed in the
adaptive mesh refinement is $0.00390625$.  The numerical constants are
given in Table \ref{tab:num_param}.  Figure \ref{fig:vortex_3d} shows
the iso-contour $\Phi = 0$ at different times. The undeformed sphere
is recovered after one cycle.

\subsection{Two Phase Flows}\label{sec:validation_two}
\subsubsection{Rising Bubble: Surface Tension Benchmark}
\begin{figure}[!h]
\centering
\begin{subfigure}[]{0.32\textwidth}
\includegraphics[width=\textwidth]{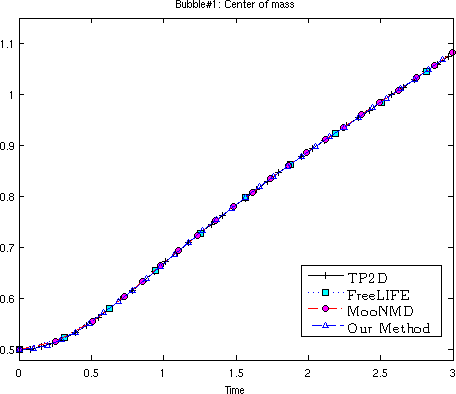}
\caption{Center of Mass, case $\#1$}
\end{subfigure}
\begin{subfigure}[]{0.32\textwidth}
\includegraphics[width=\textwidth]{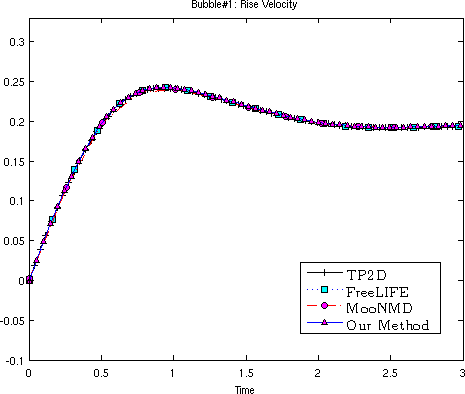}
\caption{Rise Velocity, case $\#1$}
\end{subfigure}
\begin{subfigure}[]{0.32\textwidth}
\includegraphics[width=\textwidth]{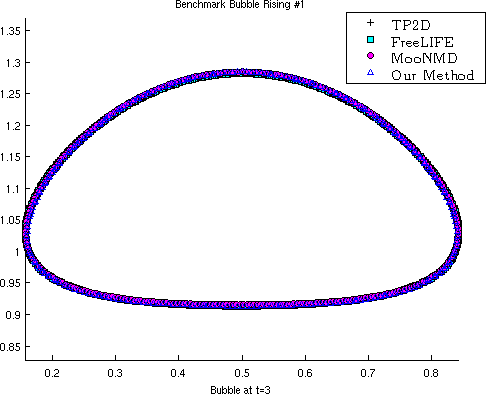}
\caption{Zero level-set at $t=3$, case $\#1$}
\end{subfigure}\\
\centering
\begin{subfigure}[]{0.32\textwidth}
\includegraphics[width=\textwidth]{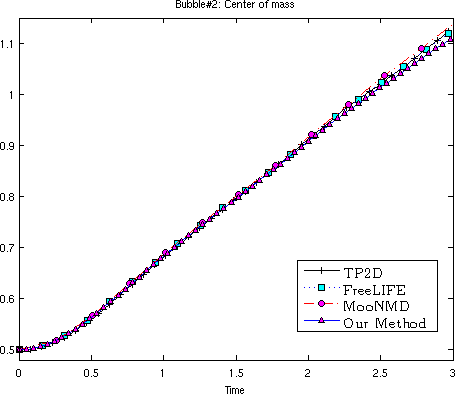}
\caption{Center of Mass, case $\#2$}
\end{subfigure}
\begin{subfigure}[]{0.32\textwidth}
\includegraphics[width=\textwidth]{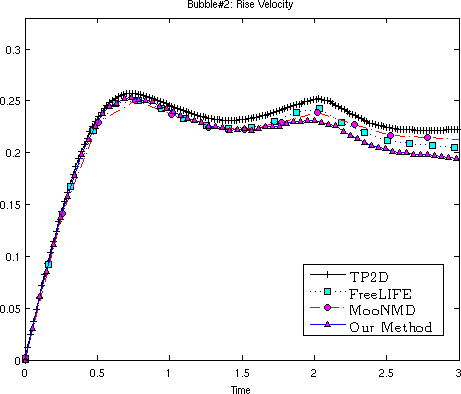}
\caption{Rise Velocity, case $\#2$}
\end{subfigure}
\begin{subfigure}[]{0.32\textwidth}
\includegraphics[width=\textwidth]{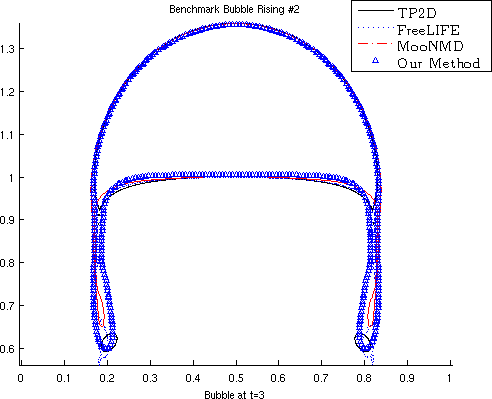}
\caption{Zero level-set at $t=3$, case $\#2$}
\end{subfigure}
\caption{Rising bubble test cases $\#1$ and $\#2$. Our simulations are
  within the range of the benchmark.}
\label{fig:ref_bubble_1}
\end{figure}

\begin{figure}[!h]
\centering
\includegraphics[scale=0.5]{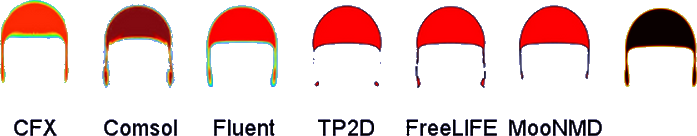}
\caption{Rising bubble case $\#2$ at $t=3$.  Different algorithms
  described in \cite{hysing2009quantitative} and our simulation
  in the rightmost panel.  The shapes of the bubbles are qualitatively
  similar.  The left figures are courtesy of S. Turek
  \cite{bubble_website}.  }
\label{fig:ref_bubble_2_2}
\end{figure}

We start the validation of the two phase flow system with the Rising
bubble benchmark problem, wee \eg \cite{hysing2009quantitative}.  The
computational domain is $\Lambda = (0,1) \CROSS (0,2)$, the initial
data $\phi_0$ is the characteristic function of a circular bubble of
radius $0.25$ centered at $(0.5,0.5)$. Two different sets of physical
constants are considered, see Table \ref{tab:benchmark}; $\sigma$ is
the surface tension coefficient and $g$ is the magnitude of the non-dimensional gravitational force.  The domains
$\Omega^+$ and $\Omega^-$ are the domains inside and outside the
bubble, respectively.  The no-slip boundary condition is imposed at
the top and bottom of the computational domain. The free slip
condition \eqref{eqn:slip} is imposed on the side walls.
\begin{table}[!ht]
\centering
\begin{tabular}{c |c | c  |  c  |  c | c | c}
    \hline
   Test case &$\rho^+$  & $\rho^-$   & $\mu^+$ &  $\mu^-$   & $g$   & $\sigma $ \\ \hline
    1        &1000      &  100         & 10      &  1.       &  0.98 & 24.5     \\  \hline
    2        &1000      &  1         & 10      &  0.1       &  0.98 & 1.96     \\  \hline
\end{tabular}
\caption{Two different sets of physical constants for the rising bubble benchmark problems.
  The first test case consists of a density/viscosity ratio of 10 with 
  a large surface tension coefficient. The density/viscosity ratio is equal to 100 in the second
  case and the surface tension coefficient is smaller.}
\label{tab:benchmark}
\end{table}
The initial subdivision is made of uniform cells of diameter
$0.03125$; the minimum mesh size allowed in the adaptive mesh
refinement is $0.0078125$. The time steps $\delta t$ are chosen
uniform according to the CFL restriction \eqref{eqn:CFL} and the
values of the numerical constants are given in Table
\ref{tab:num_param}.  We compare in Figure \ref{fig:ref_bubble_1} 
our results with those from three other
methods. We show in panel (a) the time history of the center of mass
$\bX_c := \int_{\Omega_2} \bx d\bx/\int_{\Omega_2} 1 d\bx$, in panel
(b) the rising velocity $u_c :=\int_{\Omega_2} u_y
d\bx/\int_{\Omega_2} 1 d\bx$, and in panel (c) the shape of the bubble
at $t=3$. The results of our simulations
are within the range of those given by the benchmark algorithms.
The shape of the bubble for the test case \#2 at $t=3$ is reported
in Figure \ref{fig:ref_bubble_2_2} and compared with the shapes
obtained by the other algorithms described in
\cite{hysing2009quantitative}.

\subsubsection{Buckling fluids}\label{ss:buckling}
We now test the algorithm in the context of fluid buckling, see for
instance \cite{ville2011convected, ISI:000079208200005,MR2219879}.
The test consists of letting a free-falling jet of very viscous fluid
impinge on a horizontal surface.  The diameter of the impinging jet is
\SI{0.1}{\metre} and the inflow velocity is \SI{1}{\metre\per\second}
in $\Lambda=(0,\SI{1}{\metre})^2$.  The physical parameters chosen for
the falling fluid are $\rho^+=\SI{1800}{\kilogram\per\metre\cubed}$,
$\mu^+=\SI{250}{\pascal\second}$ and
$\rho^-=\SI{1}{\kilogram\per\metre\cubed}$,
$\mu^-=\SI{2e-5}{\pascal\second}$ for the ambient fluid as in
\cite{ville2011convected}.  The no-slip boundary condition ($\bu =
\bf{0}$) is imposed at the bottom boundary, and an inflow boundary
condition is imposed where $|x-0.5|< 0.05$ at the top boundary.  The open
boundary condition is enforced on all the other boundaries, \ie $(2\mu
\nabla^S\bu - pI)\bnu = 0$.  Surface tension is neglected for this
test case ($\sigma = 0$).  The numerical constants are given in Table
\ref{tab:num_param}.

\begin{figure}[!h]
 \begin{center}
   {\includegraphics[scale=0.16]{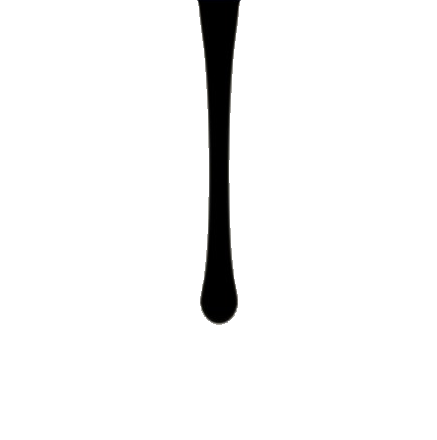}}
   {\includegraphics[scale=0.195]{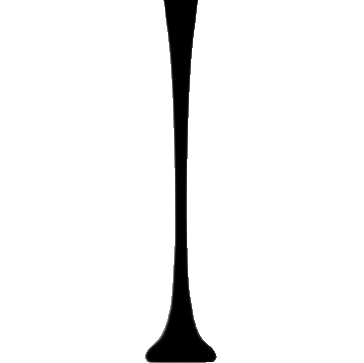}}
   {\includegraphics[scale=0.195]{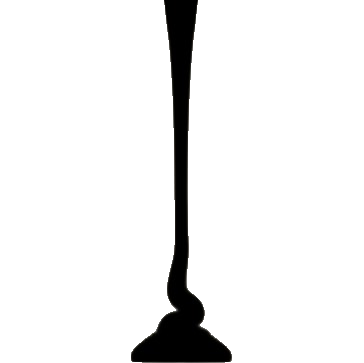}}
   {\includegraphics[scale=0.195]{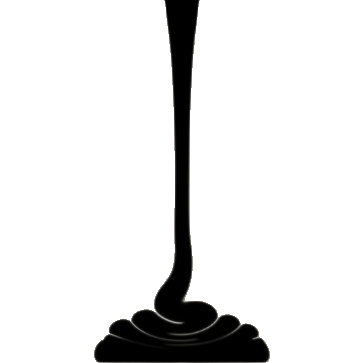}}
   {\includegraphics[scale=0.195]{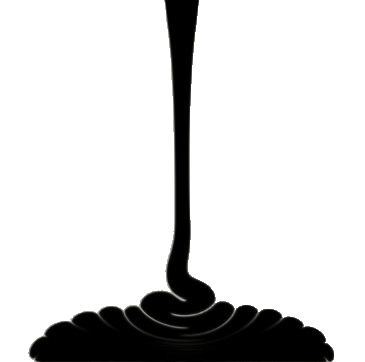}}
 \end{center}
 \caption{(From left to right) Time evolution of a jet of very
   viscous liquid inside a cavity filled with air.  Buckling occurs
   when the liquid hits the floor.}
 \label{fig:fluid_buck}
 \end{figure}

The initial subdivision is composed of cells of uniform diameter
$\SI{0.015625}{\metre}$, and the minimum cell diameter reached during
the mesh adaption process is $\SI{0.003906}{\metre}$.  The time
evolution of the fluid is shown in Figure
\ref{fig:fluid_buck}. Buckling occurs after the viscous fluid impacts
the rigid bottom plate.
 \begin{figure}[!h]
 \begin{center}
 \begin{tabular}{ccccc|c}
    {\includegraphics[scale=0.18] {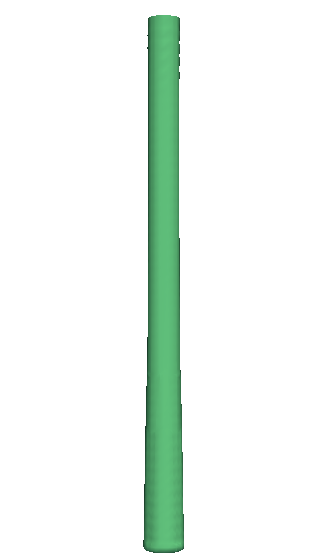}}&
    {\includegraphics[scale=0.18] {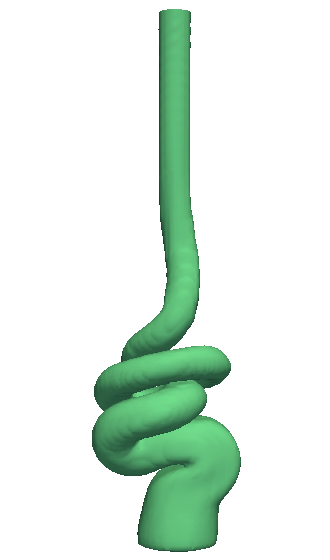}}&
    {\includegraphics[scale=0.18] {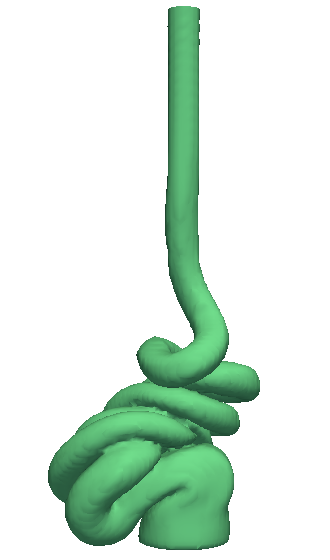}}&
    {\includegraphics[scale=0.18] {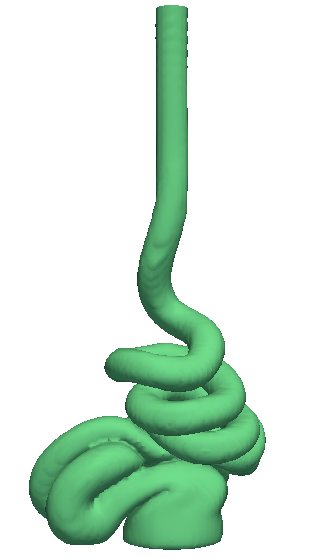}}&
    {\includegraphics[scale=0.18] {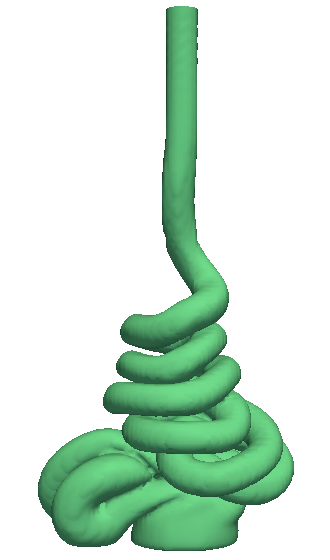}}&
    {\includegraphics[scale=0.3] {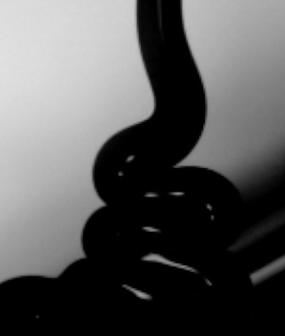}}\\
    (a) & (b) & (c) & (d) & (e) & (f)
    \end{tabular}
 \end{center}
 \caption{(a)-(e) Time evolution of a three-dimensional jet of
   Silicone oil falling in a cavity filled with air at time (a)
   $t=0.00501~s$, (b) $t=0.0209~s$, (c) $t=0.031~s$, (d)
   $t=0.04099~s$, and (e) $t=0.048998~s$. For comparison, (f) shows
   the shape of the Silicone oil jet obtained in laboratory with the
   same physical conditions.  Our numerical simulations are
   qualitatively in agreement with the physical experiment.}
 \label{fig:fluid_buck-3d}
 \end{figure}

 We also show in Figure \ref{fig:fluid_buck-3d} a three-dimensional
 simulation. The computational domain is $\Lambda =
 (0,\SI{0.008}{\metre})^3$, the initial diameter of the jet is
 $\SI{0.0004}{\metre}$ and the inflow velocity is
 $\SI{1.75}{\metre\per\second}$. The viscosity and density of the
 fluid are those of silicone oil:
 $\rho^+=\SI{960}{\kilogram\per\metre\cubed}$,
 $\mu^+=\SI{5}{\pascal\second}$. The viscosity and density of the
 ambient fluid are those of air:
 $\rho^-=\SI{1.2}{\kilogram\per\metre\cubed}$,
 $\mu^-=\SI{2e-5}{\pascal\second}$.  Surface tension is accounted for,
 $\sigma=\SI{0.021}{\newton\per\metre}$.  The no-slip boundary
 condition ($\bu = \bf{0}$) is imposed at the bottom of the box. The
 inflow boundary is the disk $\sqrt{(x-0.004)^2 + (y-0.004)^2}<
 0.0002$ at the top boundary ($z=0.008$).  Open boundary conditions
 $(2\mu \nabla^S\bu - pI)\bnu = 0$ are applied on the rest of the
 boundary.  The results are obtained with adaptive mesh refinement
 with the minimum cell diameter \SI{6.25e-5}{\metre}.  The time steps
 follow the CFL restriction \eqref{eqn:CFL} and the numerical
 constants are given in Table \ref{tab:num_param}.

\section{Numerical Simulations of Bouncing Jets}
\label{sec:bouncing_jets}
We now use our algorithm to predict the bouncing effect of jets of
Newtonian (Section \ref{subsec:newtonian_bounce_jet}) and
non-Newtonian (Section \ref{subsec:kaye_effect})fluids. In both cases,
the formation of a thin layer of air between the jet and the bulk of
the fluid is a critical ingredient to observe the bouncing effect.
The adaptive mesh refinement strategy adopted in our algorithm allows
to capture this thin layer with a reasonable number of degrees of
freedom.


\subsection{Two-dimensional Newtonian Bouncing Jets}
\label{subsec:newtonian_bounce_jet}
We start with a Newtonian fluid falling into a translating bath as in
the experiment proposed in \cite{BounceJet}.  The fluid is a silicone
oil with viscosity $\mu^+=\SI{0.25}{\pascal\second}$ and density
$\rho^+=\SI{960}{\kilogram\per\metre\cubed}$. The ambient fluid is
air: $\rho^-=\SI{1.2}{\kilogram\per\metre\cubed}$ and
$\mu^-=\SI{2e-5}{\pascal\second}$.
The computational domain is
{$\Lambda:=(0,\SI{0.04}{\centi\metre})\CROSS
  (0,\SI{0.06}{\centi\metre})$}.  The radius of the incoming jet is
\SI{0.25}{\centi\metre}; its velocity is $\bu = (0,
-\SI{5}{\centi\metre\per\second})$ on $\{(x,y)\in \partial \Lambda \ :
y=4, \ |x-0.01|<0.0025\}$. The region $\{(x,y)\in \Lambda \ : \ 0< y <
0.02\}$ is filled by the same fluid, called the ``bath'', and it moves
to the right with a horizontal velocity $V_B$ to be specified later.
Slip boundary condition \eqref{eqn:slip} are imposed at the bottom of
the cavity and the open boundary condition $(2\mu \nabla^S\bu - pI)\bnu =
0$ is applied to the rest of the boundary
$$
\{(x,y) \in \partial \Lambda \ : \ y=4, \ |x-0.01|\geq 0.0025\} \cup
\{(x,y) \in \partial \Lambda \ : \ x=0, \ y\geq 0.02\} \cup \{(x,y)
\in \partial \Lambda \ : \ x=0.04\}.
$$
The values of the numerical constants are those given in Table
\ref{tab:num_param}.  The minimum cell diameter resulting from the
adaptive mesh refinement strategy is $\SI{7.8125e-5}{\centi\metre}$.
The time steps are chosen to satisfy the CFL restriction
\eqref{eqn:CFL}.

\begin{figure}[!ht]
\centering
\begin{subfigure}[]{0.9\textwidth}
\centering
 {\includegraphics[scale=0.2]{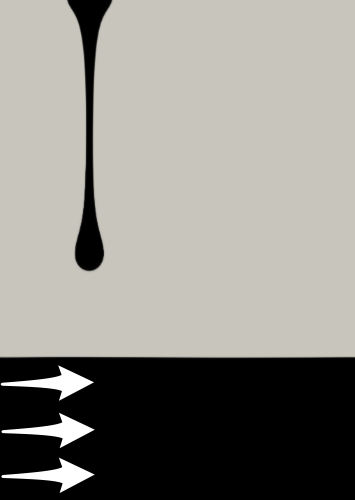}} 
{\includegraphics[scale=0.2]{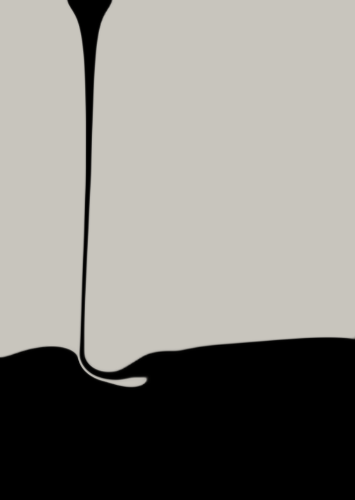}} 
{\includegraphics[scale=0.2]{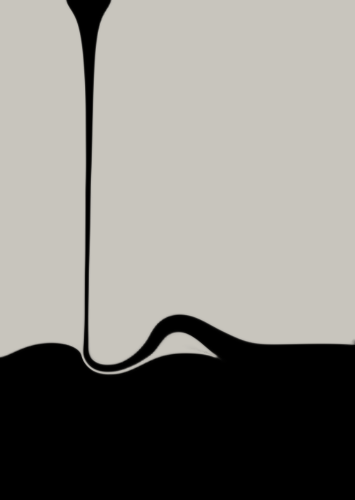}} 
{\includegraphics[scale=0.2]{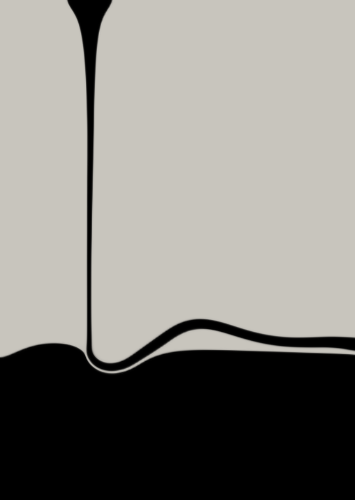}} 
\caption{}
\label{fig:bouncing_jet_a}
\end{subfigure}
\begin{subfigure}[]{0.9\textwidth}
\centering
{\includegraphics[scale=0.2]{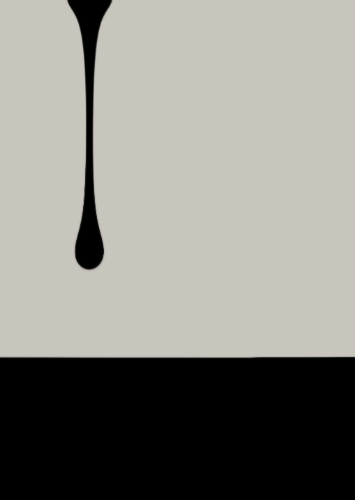}} 
{\includegraphics[scale=0.2]{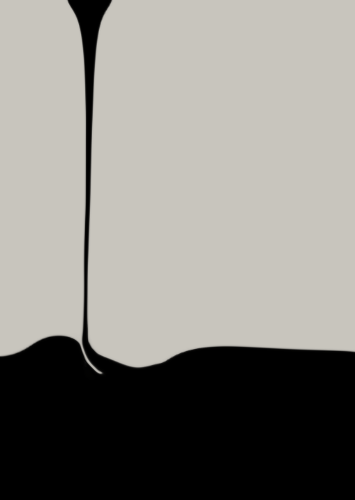}} 
{\includegraphics[scale=0.2]{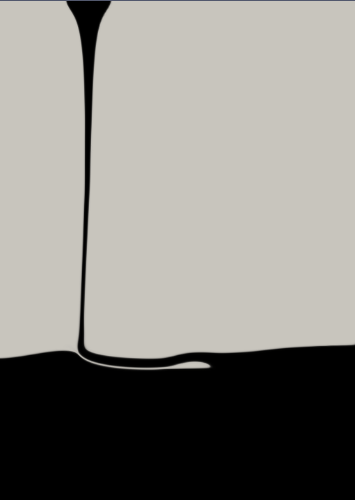}} 
{\includegraphics[scale=0.2]{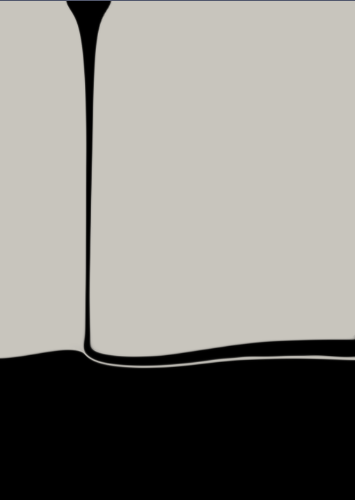}} 
\caption{}
\label{fig:bouncing_jet_b}
\end{subfigure}
\label{fig:bouncing_jet}
\begin{subfigure}[]{0.9\textwidth}
\centering
{\includegraphics[scale=0.2]{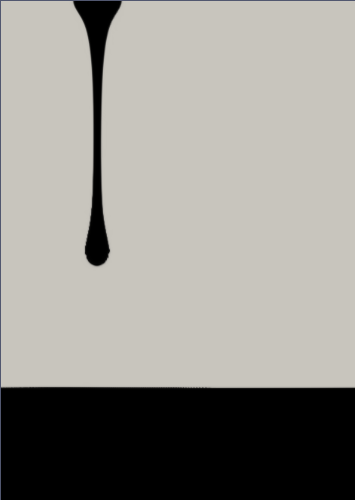}} 
{\includegraphics[scale=0.2]{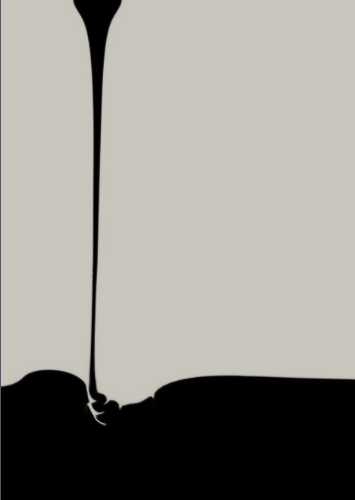}} 
{\includegraphics[scale=0.2]{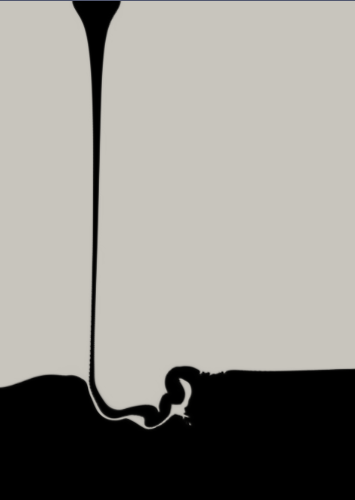}} 
{\includegraphics[scale=0.2]{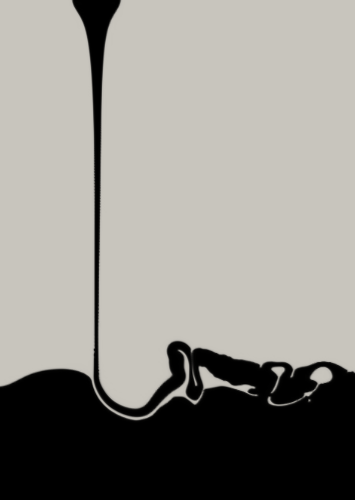}} 
\caption{}
\label{fig:bouncing_jet_c}
\end{subfigure}
\caption{Newtonian bouncing jet (from left to right). The white arrow
  indicates the translation direction of the bath.  (a) Bath velocity
  $V_B=\SI{8}{\centi\metre\per\second}$ and surface tension
  coefficient $\sigma=\SI{21}{\milli\newton\per\metre}$;  the incoming
  jet bounces away from the bath and we observe the apparition of an
  air layer between the jet and the bath.  (b) Bath velocity
  $V_B=\SI{25}{\centi\metre\per\second}$ and surface tension
  coefficient $\sigma=\SI{21}{\milli\newton\per\metre}$;  in this
  case, the bath velocity is too large and the jet slides along the
  bath surface; (c) Bath velocity
  $V_B=\SI{8}{\centi\metre\per\second}$ but without surface tension
 $\sigma=\SI{0}{\milli\newton\per\metre}$; compare with (a).  
}
\end{figure}

As already noted in \cite{BounceJet}, there is a range of velocities
for which bouncing occurs, but the jet slides along the surface of the
bath when the horizontal velocity of the bath is too high.  This is
illustrated in Figures \ref{fig:bouncing_jet_a} and
\ref{fig:bouncing_jet_b}: the jet bounces when
$V_B=\SI{8}{\centi\metre\per\second}$ (see Figure
\ref{fig:bouncing_jet_a}) but it slides when
$V_B=\SI{25}{\centi\metre\per\second}$ (see Figure
\ref{fig:bouncing_jet_b}). Note that it is necessary to include the
surface tension to keep the jet stable after impinging on the free
surface as illustrated in Figure \ref{fig:bouncing_jet_c} where a
simulation without surface tension is presented.  Observe finally that
all our numerical simulations show that a thin layer of air is formed
between the jet and the bath each time the jet
bounces.

\subsection{Kaye Effect}
\label{subsec:kaye_effect} 
\begin{figure}[!ht]
 \centering
 \begin{subfigure}[]{0.15\textwidth}
\includegraphics[width=\textwidth]{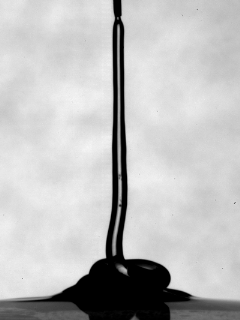}
\caption{5 ml/min}
\end{subfigure}
\begin{subfigure}[]{0.15\textwidth}
\includegraphics[width=\textwidth]{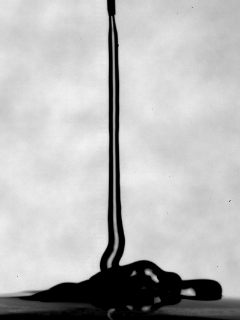}
\caption{6 ml/min}
\end{subfigure}
\begin{subfigure}[]{0.15\textwidth}
\includegraphics[width=\textwidth]{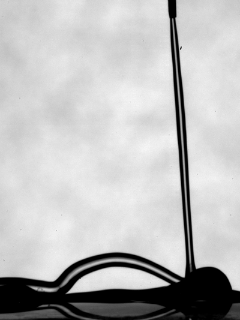}
\caption{7 ml/min}
\end{subfigure}
\begin{subfigure}[]{0.15\textwidth}
\includegraphics[width=\textwidth]{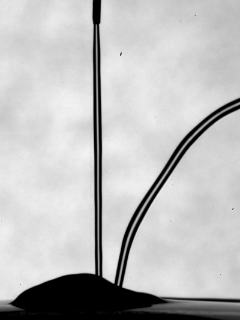}
\caption{8 ml/min}
\end{subfigure}
\caption{
  Shampoo poured at different flow rates (a)
  \SI{5}{\milli\liter\per\minute}; (b)
  \SI{6}{\milli\liter\per\minute}; (c)
  \SI{7}{\milli\liter\per\minute}; and (d)
  \SI{8}{\milli\liter\per\minute}; No bouncing is observed at
  low flow rate (a) and (b).  However, the Kaye effect is observed at
  higher flow rate (c) and (d).}
\label{fig:diff_Q}
\end{figure}

The Kaye effect is the name given to the bouncing jet phenomenon when
the bath is stationary. This effect has been observed to occur only
with non-Newtonian fluids.  It is now recognized that shear-thinning
viscosity is a critical component of the Kaye effects, see \eg
\cite{C_1976,versluis2006leaping,B_2009}.  Following
\cite{versluis2006leaping}, we adopt the
model of \cite{Cross1965417} in the rest of the paper
\begin{equation}
  \mu(\gamma) = \mu_{\infty} + \dfrac{\mu_0 - \mu_{\infty}}{1 + \left(\dfrac{\gamma}{\gamma_c}\right)^n} ,
\label{eqn:shearthinning}
\end{equation}
where $\mu_0$ is the viscosity at zero shear stress, $\mu_{\infty}$ is
the limiting viscosity for large stresses, $\gamma$ is the Frobenius norm of the rate-of-strain tensor
\[
\gamma:= \|\nabla^S \bu\|:= \left(\sum_{i=1}^d\sum_{j=1}^d \frac14  \left(\frac{\partial \bu_j }{\partial \bx_i} + \frac{\partial \bu_i}{\partial \bx_j} \right)^2\right)^{\frac12}
\] 
and $\gamma_c$, $n$ are two additional parameters.  As a benchmark, we
consider a commercial shampoo for which the shear-thinning constants
corresponding to the above model \eqref{eqn:shearthinning} have been
identified experimentally in \cite{PhysRevE.87.061001}. Recall that
viscosities cannot be measured directly but are deduced from velocity
and displacement measurements so as to match in some least-squares
sense a behavior conjectured a-priori.  The parameters $\mu_0$,
$\mu_\infty$, $\gamma_c$, $n$ corresponding to the model
\eqref{eqn:shearthinning} obtained in \cite{PhysRevE.87.061001} are
\begin{equation}
\mu_0 = \SI{5.7}{\pascal\second}, \quad \mu_{\infty} = \SI{e-3}{\pascal\second}, \quad \gamma_c = \SI{15}{\per\second}, \quad \mathrm{and} \quad  n=1.
\label{e:shear_constant}
\end{equation}
We show in Figure \ref{fig:diff_Q} snapshots of experiments done with
the shampoo poured at different flow rate. It was observed unambiguously
in \cite{PhysRevE.87.061001} that the jet slides on a lubricating air layer.

In order to reproduce qualitatively the above experiments we consider
the two-dimensional computational domain $\Lambda =
(\SI{0.496}{\meter},\SI{0.594}{\meter})\CROSS(0,\SI{0.016}{\meter})$.
The no-slip boundary condition ($\bu = \bf{0}$) is imposed at the
bottom of the computational domain and an inflow boundary condition is
imposed on the disk $\{(x,y) \ | \ |x-0.5|< 0.00021875\}$ on the top
of the box.  This correspond to a jet of radius
\SI{0.4375}{\milli\meter}. The inflow velocity is taken to be
\SI{1.75}{\meter\per\second}.  The open boundary condition $(2\mu
\nabla^S\bu - pI)\bnu = 0$ is applied on all the other boundaries.
The numerical constants used in the simulations are listed in Table
\ref{tab:num_param}.  The minimum mesh size attainable by adaptive
mesh refinement is $\SI{3.125e-5}{\meter}$, and time steps are chosen
to comply with the CFL restriction \eqref{eqn:CFL}.
The physical parameters chosen for the fluid are
$\rho^+=\SI{1020}{\kilogram\per\metre\cubed}$, the shear-thinning
viscosity constants are provided in \eqref{e:new_shear_constant}. We
take $\rho^-=\SI{1.2}{\kilogram\per\metre\cubed}$ and
$\mu^-=\SI{2e-5}{\pascal\second}$ for the air. Surface tension is
applied at the fluid/air interface with the surface tension coefficient
$\sigma=\SI{0.03}{\newton\per\metre}$.

\begin{figure}[h]
\centering
\begin{subfigure}[b]{0.3\textwidth}
\includegraphics[width=\textwidth]{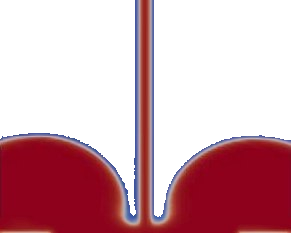}
\caption{Jet with parameters \eqref{e:shear_constant}}
\label{fig:viscosity_a}
\end{subfigure}\hfil 
\begin{subfigure}[b]{0.4\textwidth}
\includegraphics[width=\textwidth,bb=78 35  679 425,clip=]{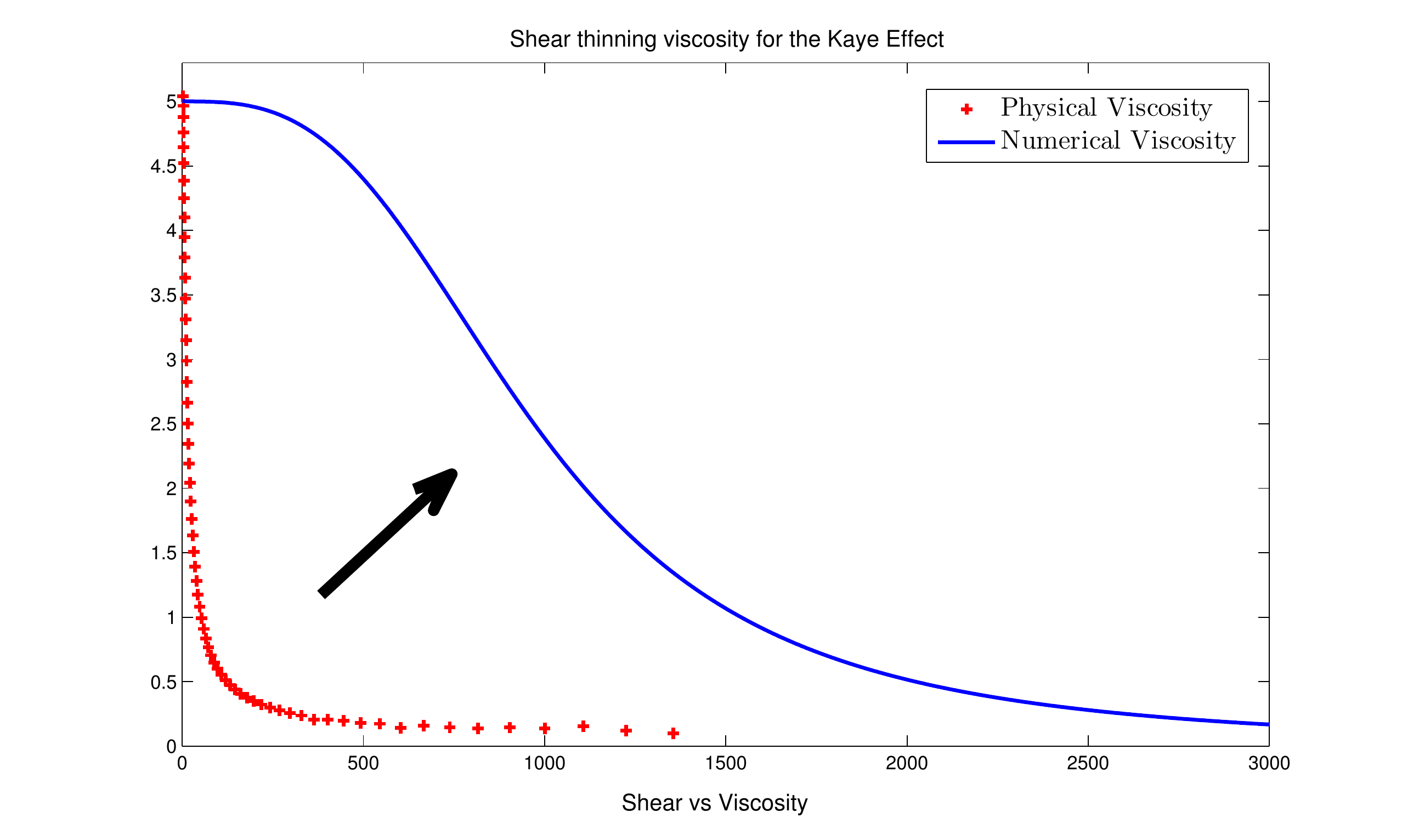}
\caption{Redefined viscosity}
\label{fig:viscosity_b}
\end{subfigure}
\caption{(Left) Falling jet of liquid described by the parameters
  given in \eqref{e:shear_constant}. The shear-thinning effect is too
  strong and the fluid becomes instantly water-like when hitting the
  bottom of the cavity.  (Right) Viscosity vs Shear for the parameters
  \eqref{e:shear_constant} (dotted line) and
  \eqref{e:new_shear_constant} (solid line).}
\label{fig:viscosity}
\end{figure}

\begin{figure}[!h]
\begin{center}
{\includegraphics[scale=0.125]{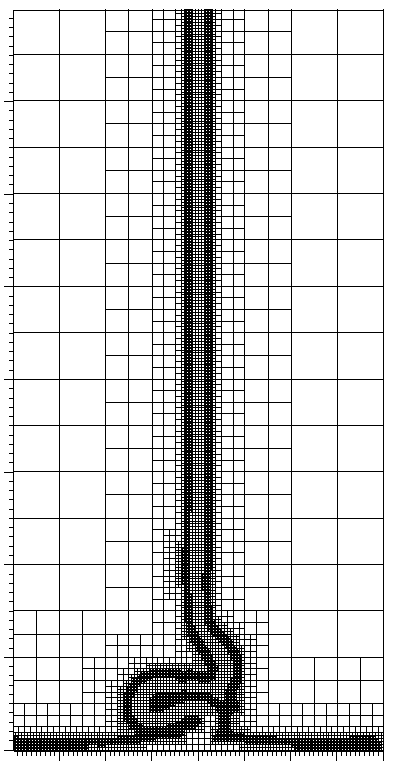}} 
{\includegraphics[scale=0.2]{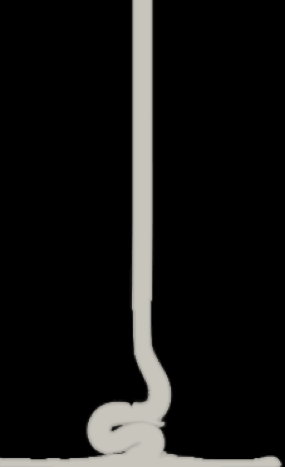}} 
{\includegraphics[scale=0.2]{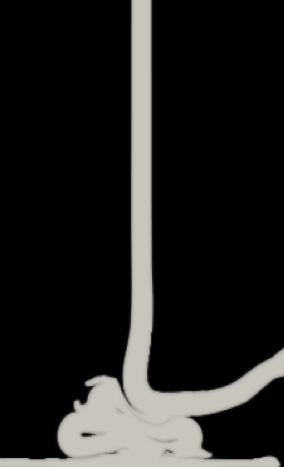}} 
{\includegraphics[scale=0.2]{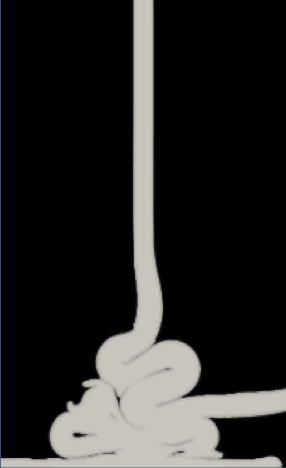}} 
{\includegraphics[scale=0.2]{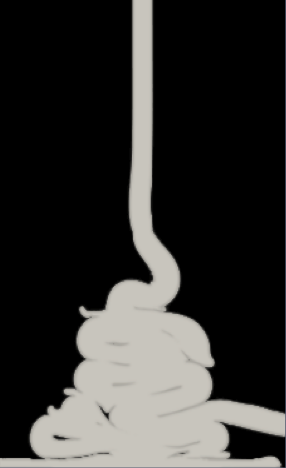}} 
{\includegraphics[scale=0.2]{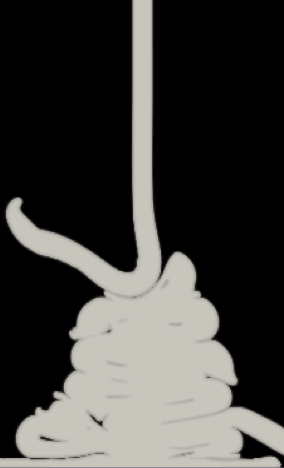}} \\
{\includegraphics[scale=0.2]{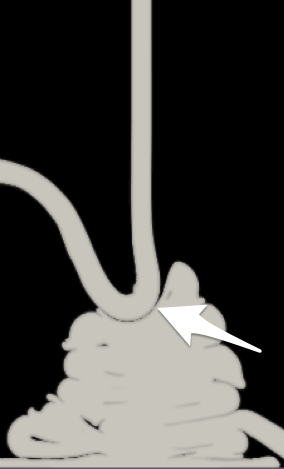}} 
{\includegraphics[scale=0.2]{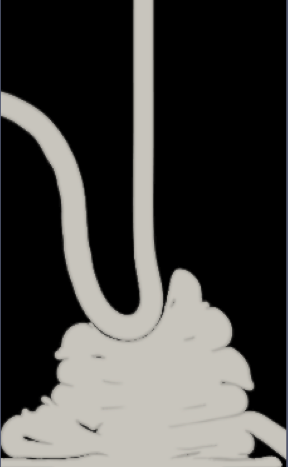}} 
{\includegraphics[scale=0.2]{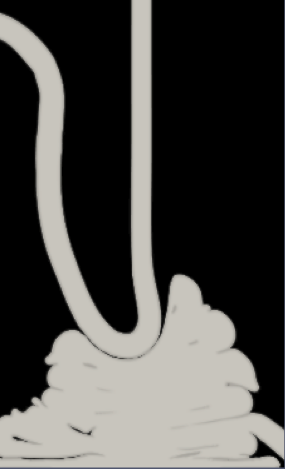}} 
{\includegraphics[scale=0.2]{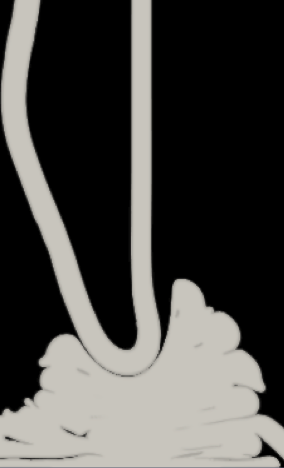}} 
{\includegraphics[scale=0.2]{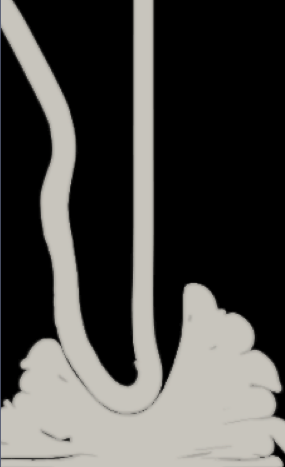}} 
{\includegraphics[scale=0.125]{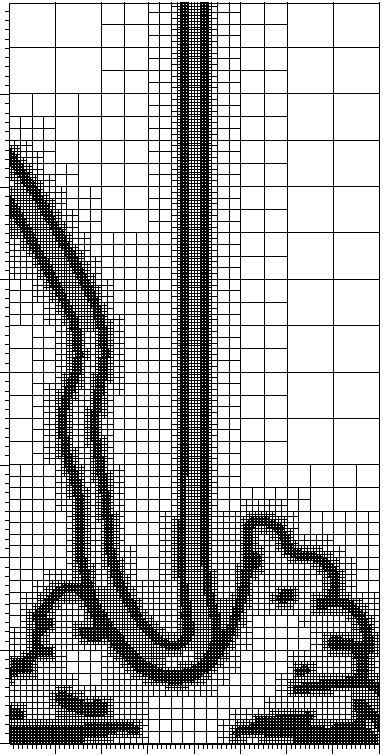}} 
\end{center}
\caption{(Left to Right, Top to Bottom) Numerical simulation of the
  Kaye effect with adaptive meshes (from left to right and top to
  bottom). The viscosity parameters are given in
  \eqref{e:new_shear_constant}.  The first and last frames illustrate
  the adaptive subdivision generated by the adaptive strategy.  An air
  layer appears in the last frame in the first row and is
  fully developed in the first frame of the second row. The apparition
  of the air layer coincides with the beginning of the bouncing
  phenomenon.}
\label{fig:num_kaye_effect}
\end{figure}

After many numerical experiments it turned out that the above physical
parameters did not give any Kaye effect.  Our interpretation is that
the shear-thinning effect given by the set of parameters
\eqref{e:shear_constant} is too strong; with these coefficients the
fluid instantly becomes water-like when hitting the bottom of the
cavity, see Figure \ref{fig:viscosity_a}.  Several explanations for
this mismatch are plausible: (i) as mentioned above the parameters
\eqref{e:shear_constant} are measured under the assumption that the
shear-thinning law follows the Cross model; (ii) the shear is
over-predicted by our algorithm; (iii) the air layer width for this
range of parameters is too thin to be captured by the algorithm; (iv)
a fundamental component is missing in our mathematical model.  These
observations have lead us to consider a different set of
shear-thinning parameters requiring a larger shear for a notable
reduction of the viscosity and a smoother transition from the maximum
to the minimum values of the viscosity, see Figure
\ref{fig:viscosity_b}.  The parameters that
we now consider are the following:
\begin{equation}
  \mu_0 = \SI{5.7}{\pascal\second}, \quad \mu_{\infty} = \SI{e-3}{\pascal\second}, 
\quad \gamma_c = \SI{970}{\per\second} \quad \mathrm{and} \quad  n=3.
\label{e:new_shear_constant}
\end{equation}
This set of parameter produces the Kaye effect as shown in Figure
\ref{fig:num_kaye_effect}.
%
Here again we notice that there is a very thin layer of air between
the bouncing jet and the rest of the fluid, thereby adding to the
large body of evidence pointing at the importance of air layers in
bouncing jets.  Let us finally mention that these computations show
also that mesh adaptivity is critical to reproduce numerically the
Kaye effect.

\section*{Acknowledgments}
The experiments reported in the paper were done at the High Speed
Fluid Imaging Laboratory of S. Thoroddsen at KAUST during two visits
of S.L. The authors would like to express their gratitude to
S. Thoroddsen and E. Li for their help in conducting these
experiments.  W. Bangerth availability and his constant help in the
implementation of our algorithms with the deal.II library is
acknowledged.  The authors acknowledge the Texas A\&M University
Brazos HPC cluster that contributed to the research reported here.

\bibliographystyle{abbrvnat}







\end{document}